%% file: IntegrationPaperMaster.tex
\documentclass[11pt]{amsart}

\usepackage{amssymb, amsmath, amscd, amsthm, color, epsfig}

\usepackage[all]{xy}          
\xyoption{dvips}              

\newcommand{\R}{{\mathbb R}}

\newcommand{\cd}[1]{\mathcal{D}^{c}_{#1}}

\newcommand{\td}[1]{\mathcal{D}_{#1}}

\newcommand{\Td}{\mathcal{D}}
\newcommand{\G}[2]{\Gamma_{#1,#2}}

\newcommand{\St}[2]{\mathcal{S}_{#1,#2}}
\newcommand{\Ss}{\mathcal{S}_{S}}
\newcommand{\ST}{\mathcal{S}_{T}}
\newcommand{\SU}{\mathcal{S}_{U}}
\newcommand{\SI}{\mathcal{S}_{I}}
\newcommand{\SH}{\mathcal{S}_{H}}
\newcommand{\SX}{\mathcal{S}_{X}}

\newcommand{\Sd}{\mathcal{S}_{disc}}

\newcommand{\K}{{\mathcal{K}}}
\newcommand{\W}{{\mathcal{W}}}
\newcommand{\V}{{\mathcal{V}}}

\newcommand{\FM}[2]{F[#1, #2]}

\theoremstyle{plain}
\newtheorem{thm}{Theorem}[section]
\newtheorem{prop}[thm]{Proposition}
\newtheorem{lemma}[thm]{Lemma}
\newtheorem{cor}[thm]{Corollary}

\theoremstyle{definition}
\newtheorem{definition}[thm]{Definition}

\theoremstyle{remark}
\newtheorem*{rem}{Remark}
\newtheorem*{rems}{Remarks}

\newcommand{\refT}[1]{Theorem~\ref{T:#1}}
\newcommand{\refC}[1]{Corollary~\ref{C:#1}}
\newcommand{\refP}[1]{Proposition~\ref{P:#1}}
\newcommand{\refD}[1]{Definition~\ref{D:#1}}
\newcommand{\refL}[1]{Lemma~\ref{L:#1}}

\begin{document}

\title{A survey of Bott-Taubes integration}
\author{Ismar Voli\'c}
\address{Department of Mathematics, University of Virginia, 
Charlottesville, VA}
\email{ismar@virginia.edu}
\urladdr{http://www.people.virginia.edu/\~{}iv2n}
\subjclass{Primary: 57M27; Secondary: 81Q30, 57R40}
\keywords{configuration space integrals, knots, spaces of knots,
 Feynman diagrams, Chern-Simons theory, finite type invariants}

\begin{abstract}
It is well known that certain combinations of configuration space integrals defined by Bott and Taubes \cite{BT} produce cohomology classes of spaces of knots.  The literature surrounding this important fact, however, is somewhat incomplete and lacking in detail.  The aim of this paper is to fill in the gaps as well as summarize the importance of these integrals.
\end{abstract}

\maketitle

{\tableofcontents}

  \section{Introduction}\label{S:Intro}

Let $\K_m$ be the space of smooth embeddings of $S^1$ in $\R^m$ and set $\K_3=\K$ to simplify notation.
 There is a well-known way to produce cohomology classes on 
 $\K_{m}$ coming from the the perturbative  
 Chern-Simons theory.  In case of classical knots $\K$, the construction yields special classes 
 in zeroth cohomology, finite type 
 knot invariants, which have been studied extensively over the past 
 ten years by both topologists and physicists.  The physics behind the 
 construction will not be 
 discussed here, but a good start would be the informative survey 
  by Labastida \cite{Lab}.

   At the heart of the theory are 
  configuration space integrals which first appeared some fifteen years ago in the work of Guadagnini, Martellini, and Mintchev \cite{GMM} and Bar-Natan \cite{BN3}.  The version we are concerned with here is more recent and is due to Bott and Taubes \cite{BT}.  The idea is that certain chord diagrams may be used as prescriptions for obtaining a form on a configuration space of points in $\R_m$, some of which are required to lie on a given knot, and then pushing this form forward to $\K_m$.  Adding the result over all diagrams gives a closed form on $\K_m$.  For nice overviews of the constructions in case of $\K$, see \cite{Bott1, Bott2, L2}.
  
  The most complete proofs of this fact, which has been used very effectively in recent years, are due to Altschuler and Freidel \cite{Alt}.  Among other things, they were the first to prove that this approach produces a universal finite type invariant (see \S\ref{S:Universal}).  However, their arguments are 
inspired
from physics and may not be transparent to an algebraic topologist.  Same is true of the important work of Poirier \cite{Poir, Poir2} who extended Bott-Taubes integrals to links and tangles (also see related work of Yang \cite{Y}). 

On the other hand, the more topological arguments, which can be found in the original Bott-Taubes work as well as in D. Thurston's undergraduate thesis \cite{Th}, are somewhat incomplete even if most of the necessary ideas are there.  Further, Thurston's work was unfortunately never published.  This paper owes much to Bott-Taubes and Thurston, and most of the arguments used here come directly from, or are inspired by, the ideas they give in \cite{BT} and \cite{Th}. 

The focus of Thurston's work is on $\K$, but the generalization to $\K_m$, originally due to Cattaneo, Cotta-Ramusino, and Longoni 
\cite{Catt}, is straightforward.  Depending on $m$, however, one now obtains cohomology classes in different dimensions.  
The goal of this paper is to give detailed proofs of these facts.  The statements we are after are \refT{Thurston}, \refT{Universality}, and a weak version of \refT{CCL}, stated here loosely together as:

\begin{thm}\label{T:IntroMain}
Bott-Taubes configuration space integrals combine to yield nontrivial cohomology classes of $\K_m$.  For $\K$, they represent a universal finite type knot invariant.
\end{thm}
The bulk of the paper, \S\ref{S:Vanishing}, is devoted to proving 
the first part of this theorem.   The proof is essentially Stokes' Theorem, where integration along various types 
of codimension one faces of the compactified configuration space has to be considered.  The goal is to show that these integrals either vanish of cancel to give a 
closed form on $\K_m$.  We do all this for $\K$ and then note in \S\ref{S:HigherCodimension} that one only needs to introduce minor changes throughout for $\K_m$.  The second part of \refT{IntroMain} is addressed in \S\ref{S:Universal}, where we also review the basics of finite type
theory.  The material in \S\ref{S:Vanishing} and \S\ref{S:Universal} will hopefully fill 
the gap in literature mentioned above and provide a more topological treatment.

Along the way, we also introduce the algebra of trivalent diagrams in \S\ref{S:Trivalent}, motivate the definition of Bott-Taubes integrals in \S\ref{S:LinkingNumber}, and then take a necessary digression on Fulton-MacPherson compactification of configuration spaces in \S\ref{S:F-M}.  The integration itself is described in 
\S\ref{S:IntegralsCohomology}.

\vskip 4pt
\noindent
Bott-Taubes configuration space integrals have played an important role in some recent developments in knot theory and beyond.  For example, Cattaneo et. al. constructed a double complex from the set of chord diagrams and used this to show that, for any $k$, $\K_m$ has nontrivial cohomology in degrees greater than $k$ for $m>3$.  Further, their double complex has has been related to the cohomology spectral sequence set up by Vassiliev \cite{Vas} and converging to $\K_m$ \cite{Tur}.  This was used in \cite{Vo2} for showing the collapse of Vassiliev's spectral sequence along a certain line.  In general, the hope is that Bott-Taubes integrals produce all cohomology classes of $\K_m$.  

This spectral sequence result uses an extension of these integrals to ``punctured knots" making up the stages of certain towers of spaces approximating $\K_m$ (in the sense of calculus of functors \cite{Vo2}).  This extension was also used to place finite type knot invariants in a more homotopy-theoretic framework \cite{Vo}.  We feel this work would benefit from a firmer grounding which we attempt to provide here.

Another active area of investigation is the relationship between Bott-Taubes integrals and the \emph{Kontsevich Integral} \cite{Kont}, which is essentially the only other known way of producing a universal finite type invariant.  Poirier characterized the relationship between the two approaches and reduced the question of their equivalence to the computation of a certain anomalous term (see discussion following \refP{anomalous}).  It is still an open and interesting question if this anomaly vanishes.  

Analogous comparison occurs between various ways of producing invariants of rational homology 3-spheres and knots in rational homology 3-spheres.  One approach is through a generalization of Bott-Taubes integration, as done by Bott and Cattaneo \cite{BC, BC2}, while the other ones are based on the Kontsevich Integral and come in many related variants \cite{LMO, BGRT, KT}.  The relationship between the two is still not well understood.

\subsection{Acknowledgements}  The author would like to thank Riccardo Longoni, Greg Arone, Pascal Lambrechts, and Dev Sinha for helpful
 comments and suggestions.


\section{Trivalent diagrams}\label{S:Trivalent}

Before we talk about configuration space integrals, we 
introduce a class of diagrams which turns 
out to nicely keep track of the combinatorics associated to those 
integrals.  These diagrams (sometimes called Feynman) are 
frequent in Chern-Simons theory.

\begin{definition}\label{D:TrivalentDiagrams}  A \emph{trivalent 
diagram of degree $n$} is a connected graph composed of an oriented interval 
and some number of line segments connecting $2n$ vertices of 
two types:
\begin{itemize}
\item
\emph{interval vertices}, constrained to lie on the interval, 
from which only one line segment emanates, 
and
\item
\emph{free vertices}, or those not constrained to lie on the interval, 
from which exactly three line segments emanate.
\end{itemize}
If a line segment connects two interval vertices, it is called a 
\emph{chord}; otherwise it is called an \emph{edge}.  The vertices are labeled $1, ..., 2n$, and each chord and edge is also oriented.
\end{definition}

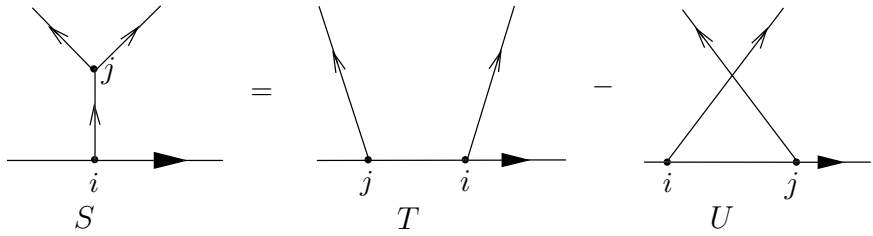
\begin{figure}[h]
\begin{center}
\input{LabeledSTU3.pstex_t}
\caption{$STU$ relation.  The three diagrams here, as well as in the $IHX$ relation below, are identical outside the 
pictured portions.}\label{F:STU}
\end{center}
\end{figure}

Let $TD_{n}$ be the set of all 
trivalent diagrams of degree $n$.  Let $STU$ be the relation 
as in Figure \ref{F:STU}.

\begin{definition}\label{D:TrivalentSpace} Let $\Td_n$ be the real vector space 
generated by $TD_{n}$, modulo the $STU$ relation.
\end{definition}


We next state a result, due to Bar-Natan \cite{BN} (see Theorem 6), which will be 
useful later and follows from the $STU$ relation. 

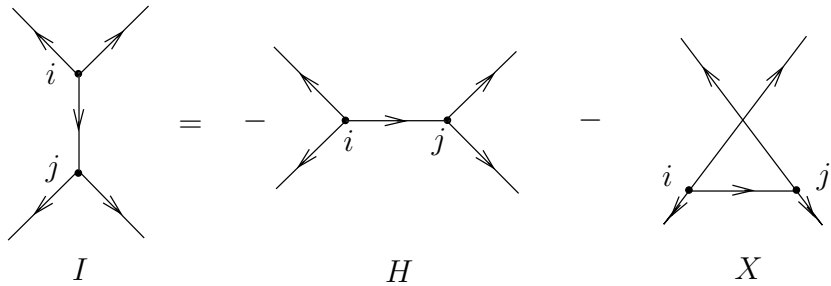
\begin{figure}[h]
\begin{center}
\input{LabeledIHX2.pstex_t}
\caption{$IHX$ relation}\label{Fi:IHXPicture}
\end{center}
\end{figure}

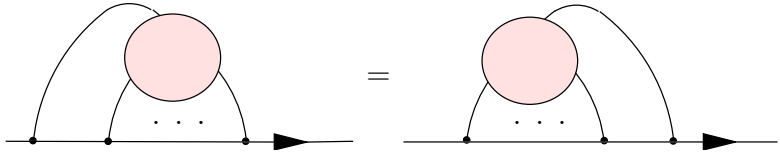
\begin{figure}[h]
\begin{center}
\input{Closure.pstex_t}
\caption{Closure relation}\label{Fi:ClosurePicture}
\end{center}
\end{figure}

\begin{prop}\label{P:IHX}  The identities given in Figures \ref{Fi:IHXPicture} and \ref{Fi:ClosurePicture} also hold in $\Td_{n}$.  The second figure is meant to indicate 
that if $D_{1}$ and $D_{2}$ ``close'' by identification 
of the endpoints of the 
interval to the same circular 
diagram considered up to orientation-preserving diffeomorphism of the 
circle, then $D_{1}=D_{2}$.
\end{prop}


\noindent


Labeling a diagram and orienting its chords and edges in different 
ways produces a 
potentially large number of non-isomorphic diagrams.  We reduce this 
number by requiring that the chords and edges are always oriented in 
such a way that they point from the vertex with lower label to the one 
with the higher label.  

Let $\Td=\oplus_{n>0}\Td_{n}$.
There is a
pairing on $\Td$,
given by the operation of connected sum.  Two diagrams are multiplied
 by continuing the 
interval of one into the interval of the other (with compatible 
orientation).  Even though there are two ways to do this, the last relation 
in \refP{IHX} ensures that the resulting 
diagrams are the same, so that the operation is commutative.  The empty 
diagram serves as the identity.  This 
operation is in correspondence 
with the operation of connected sum of knots.

\begin{definition} A trivalent diagram which can be written as a product of 
two non-empty diagrams is \emph{reducible}.  Otherwise, the diagram is 
 \emph{prime}.
\end{definition}

As it turns out, 
$\Td$ also admits a coproduct, whose precise definition can be found 
in \cite{BN}.  Bar-Natan also proves
\begin{thm}[\cite{BN}, Theorem 7]\label{T:HopfAlgebra}  $\Td$ is a 
commutative and co-commutative Hopf algebra.
\end{thm}
Now define the space of \emph{weight systems} $\W$ as the dual of 
$\Td$ (and $\W_{n}$ as its \emph{degree $n$} part).  It also has the 
structure of a Hopf algebra. To understand $\W$, it therefore suffices 
to understand its primitive elements since they generate the whole 
algebra.  To that end, we have the following useful statement, also 
due to Bar-Natan:

\begin{prop}\label{P:PrimitiveWeightSystems}  
A weight system $W\in\W$ is primitive if and only if it vanishes on 
reducible diagrams.
\end{prop}
It thus suffices to consider weight systems only on prime 
diagrams since all others are obtained from these.  This will be important in \S\ref{S:HiddenFaces} 
and \S\ref{S:PrincipalFaces}.

\section{Bott-Taubes Integrals}\label{S:B-TIntegrals}

Configuration spaces arise naturally in many 
embedding questions.  The reason is that any embedding  
gives rise to a map of configuration spaces by evaluation of the 
embedding on some number of points.  This is illustrated below by the 
simple example of the linking number of two classical knots.  Bott-Taubes 
integrals can be thought of as 
generalizations of the Gauss 
integral which computes this integer invariant of 2-component links.  We also use the next section to set the notation and make some basic definitions.


\subsection{The linking number}\label{S:LinkingNumber}

Recall that, given a space $X$, the 
\emph{configuration 
space} $F(k,X)$ is the subspace of $X^{k}$ consisting of $k$ distinct,
 ordered points in 
$X$:
$$
F(k,X)=\{(x_{1},\ldots,x_{k})\in X^{k}, \ \ x_{i}\neq x_{j} \ 
\text{for all } i\neq j\}.
$$
Now
let $F(1,1;S^{1},S^{1})$ be the configuration space of 
one point on each of two disjoint circles, which is
just the torus $T^{2}$, and let $K_{1}\amalg K_{2}$ be an embedding 
of $S^{1}\amalg S^{1}$ in 
$\R^{3}$.  The evaluation map gives a composition
\begin{gather}\label{E:Linking}
 F(1,1;S^{1},S^{1})  \longrightarrow  F(2,\R^{3}) 
\longrightarrow  S^{2}\\
 (x_{1},x_{2})  \longmapsto  (K_{1}(x_{1}),K_{2}(x_{2}))  
\longmapsto 
\frac{K_{2}(x_{2})-K_{1}(x_{1})}{|K_{2}(x_{2})-K_{1}(x_{1})|}, \notag
\end{gather}
which we denote by by $h_{12}$.

Now let $\omega_{12}$ be the standard rotation-invariant unit volume 
2-form 
on $S^{2}$,
\begin{equation}\label{E:VolumeForm}
\omega_{12} = \frac{x\, dydz-y\, dxdz+z\, dxdy}{4\pi 
(x^{2}+y^{2}+z^{2})^{3/2}}.
\end{equation}
The linking number is then defined as the integral of the pullback 
of $\omega_{12}$ by $h_{12}$:

$$L(K_{1},K_{2})=\int\limits_{T^{2}}h_{12}^{*}\omega_{12}.$$
If the same 
procedure is applied to a knot, namely if we start with the 
configuration space of two points on a single circle, this 
$F(2,S^{1})$ 
is an open cylinder so that the integration may not be defined.  This can be remedied somewhat by adding the boundary to 
the cylinder.  Geometrically, 
this corresponds to allowing $x_{1}$ and 
$x_{2}$ to ``collide'' (or ``come together,'' or ``degenerate,'' as 
we will interchangeably say).  We can now use Stokes' Theorem to say
 that the Gauss integral will yield an invariant if the integral 
over the boundary of the cylinder
of the restriction of $h_{12}^{*}\omega_{12}$ to that boundary
vanishes.
The restriction, which will be refered to 
throughout as \emph{the tangential form}, is the pullback of 
$\omega_{12}$ via the extension of $h_{12}$ to the boundary.  We 
denote this 
extension by $\tau$, which is readily seen to be given by
\begin{equation}\label{E:Tangentialmap}
\tau=\frac{K'(x_{1})}{|K'(x_{1})|}.
\end{equation}
However, the integral of $\tau^{*}\omega_{12}$ does not vanish so that 
the Gauss integral does not produce an invariant.
 One can 
now look for a ``correction term'' which would exactly cancel the
contribution of the boundary integral. 
In this case, this term turns out to be 
given by the framing number of the knot \cite{Mosk}.

This search for the ``correction term'' applies in the more general 
setup 
introduced by Bott and Taubes.  Starting with a 
configuration space of any even number of points on the circle, they 
examine the various products of maps $h_{ij}$ (analogous to $h_{12}$ 
from 
above).  To get nontrivial pullbacks, the configuration spaces first 
have to be compactified.  However, in the case of more than two 
configuration points, the simplest compactification, which puts 
the fat diagonal back in, does not yield a boundary to which the 
pullbacks extend (the fat diagonal in $X^k$ is the set $\{(x_1, ..., x_k)\in X^k\colon \exists i \neq j \text { with } x_i=x_j\}$).  What is needed, it turns out, is the more 
sophisticated compactification described in the next section.  After 
this compactification is introduced, we will look at the Bott-Taubes 
method in more detail.

\subsection{Fulton-MacPherson compactification of configuration 
spaces}\label{S:F-M}

The compactification of configuration spaces we use was first defined 
in the 
setting of algebraic geometry by Fulton and MacPherson \cite{FM}.  
It was then defined
 for manifolds by Axelrod and Singer \cite{AS} who used it
in the context of Chern-Simons theory.
We recall 
here the main features of the Axelrod-Singer version and use their constructions and definitions.  Our notation, however, follows Sinha \cite{Dev1}, whose alternative definition for the compactification does not use blowups.  This makes his constructions similar to the one given by Kontsevich in 
\cite{Kont2} and perhaps more approachable for some purposes.





%

%
Let $N$ be a smooth manifold.  Another 
way to think of $F(k,N)$ is as an ordered product of $k$ copies 
on $N$ with all diagonals 
removed.  The diagonals may be indexed by the sets $S\subseteq\{1, 
\ldots, k\}$ of cardinality at least 2.  So let $\{x_{1}, x_{2}, 
\ldots, x_{k}\}$ be a point in $F(k,N)$ and let $\Delta_{S}$ denote 
the diagonal in $N^{k}$ where $x_{i}=x_{j}$ for all $i,j\in S$.
Finally let $Bl(N^{k}, \Delta_{S})$ be the \emph{blowup} of $N^{k}$ along 
$\Delta_{S}$, namely a replacement of $\Delta_{S}$ in $N^{k}$ by 
its unit normal bundle.  Since the interior of $Bl(N^{k}, \Delta_{S})$ is $F(|S|, N)$, there are natural projections of $F(k,N)$ to $Bl(N^{k}, \Delta_{S})$ for all $S$. Along with the inclusion of $F(k,N)$ in $N^{k}$, one then has an embedding of $F(k,N)$ in
$$
N^{k}\times\prod_{\substack{S\subseteq\{1,\ldots,k\} \\ |S|\geq 2}}
Bl(N^{k},\Delta_{S}).
$$

\begin{definition}\label{D:FM}  The \emph{Fulton-MacPherson 
compactification} of $F(k,N)$, denoted by 
    $F[k,N]$, is the 
closure of 
$F(k,N)$ in 
the above product.
\end{definition}
When configuration points come together, the 
directions of approach are kept track of in $F[k,N]$.  This is because 
the normal bundle records the tangent vectors of paths approaching a 
diagonal, but up to translation on the diagonal itself.  In addition 
to translation, tangent vectors are also taken up to scaling because 
the normal bundle is replaced by its sphere bundle.  This geometric 
point of view
 will be important in understanding 
the coordinates on $F[k,N]$ to be defined in 
\S\ref{S:Coordinates&Forms}.

We now list a few important properties of this compactification.  
Detailed proofs can be found in 
\cite{AS, Dev1}.  For a more succinct treatment, see 
\cite{Dev}.

\begin{itemize}
\item The inclusion of $F(k,N)$ into $F[k,N]$ induces a homotopy 
equivalence.
\item $F[k,N]$ is a smooth $k$-dimensional manifold with corners,
and it is compact if $N$ is compact.
\item  Any embedding of $M$ in $N$ induces an embedding of $F[k,M]$ 
in $F[k,N]$.
\end{itemize}

Since we will be integrating certain forms along the various codimension 
1 faces of $F[k,N]$, it is worth discussing 
the stratification of this space in a little more detail.  More 
explicit description can be found in \cite{AS}.

Intuitively, a stratum (face, screen) corresponds to a subset 
$A$ of the $k$
configuration points degenerating.  We denote this stratum by 
$\mathcal{S}_{A}$ and set $|A|=a$.  To simplify notation, also 
assume the colliding points $x_{i}$ are indexed by $1, 2, \ldots, a$.
Either all $x_{i}$ 
degenerate to the same point, or disjoint subsets of 
$(x_{1}, \ldots, x_{a})$ degenerate to different points.  Further, some
points may degenerate faster than the others.
The various ways of collisions of configuration 
points can be efficiently kept track of by \emph{nested subsets} of 
$\{1,\ldots,a\}$.  These are sets of subsets of $\{1,\ldots,a\}$, 
where each of 
the subsets contains at least two elements and each two subsets are 
either disjoint or one is contained in the other.  Every stratum can 
be characterized in this way.  So suppose that
$\mathcal{S}_{A}$ is described by $i$ nested subsets $\{A_{1},A_{2}, 
\ldots, A_{i}\}$.  Axelrod and Singer then show that
the codimension of $\mathcal{S}_{A}$ in $F[k,N]$ is $i$.  This is a direct consequence of an explicit description of coordinates on $F[k,N]$ (see \S\ref{S:Coordinates&Forms}).
It follows that the codimension one strata consist of limits of 
sequences of points degenerating to the same point, at the same 
rate.

\begin{rem}
Given the representation of the strata as nested subsets, 
it is easy to see that $\mathcal{S}_{A}$ will be contained in
$\mathcal{S}_{A'}$ as one of its faces if the first set of nested 
subsets is contained in the second.  This can be used to graphically 
depict the relationships between the various strata of $F[k,N]$.  
If $N$ is the interval, the pattern that emerges is that of Stasheff associahedra, which also 
appear in some combinatorial descriptions of finite type knot 
invariants.  Alternatively, Sinha introduces a category of 
trees in \cite{Dev1} to study the stratification.
\end{rem}

\subsection{Integrals and knot invariants}\label{S:IntegralsCohomology}

To simplify the combinatorics to come, now let $\K$ be the space of maps of the unit interval $I$ in $\R^3$ (or $S^{3}$; it will be clear form the context which we mean) which are embeddings except at endpoints and send the endpoints to the same point with the same derivative.  This is clearly the same as the ordinary space of knots. To produce 
Gauss-like integrals that yield 0-forms, or knot invariants, Bott and Taubes first
consider $F(4,S^1)$.  

Given an embedding $K\in\K$, there are maps
$$
h_{ij}\colon \FM{4}{I}\times\K \longrightarrow S^{2}, \ \ \ \ 
i,j\in\{1,2,3,4\}, \ i\neq j,
$$
given by
$$
(x_{1}, x_{2}, x_{3}, x_{4})\longmapsto 
\frac{K(x_{j})-K(x_{i})}{|K(x_{j})-K(x_{i})|}.
$$
Since $\FM{4}{I}$ is four-dimensional, take a product of two of 
the 
maps $h_{ij}$ so that the target, $S^{2}\times S^{2}$ is also
four-dimensional.  Since Bott and Taubes do not consider 
framed knots, the only interesting choice (see \S\ref{S:Universal} for an
explanation) turns out to be
$
h=h_{13}\times h_{24}.
$
Let $
\alpha=h^{*}(\omega_{13}\omega_{24})
$
be the 4-form obtained by pulling back the product of unit volume forms
$\omega_{13}$ and $\omega_{24}$ from 
$S^{2}\times S^{2}$ via $h$.
Probably the most important feature of 
the Fulton-MacPherson compactification for us is that $\alpha$ extends 
smoothly to the boundary of $\FM{k}{I}$ for any $k$, as we will show in 
\S\ref{S:Coordinates&Forms}.

Let $\pi$ be the projection of $\FM{4}{I}\times \K$ to $\K$.  Since 
this is a trivial bundle over $\K$ with 
4-dimensional compact fiber $\FM{4}{I}$, it makes sense to integrate 
$\alpha$ along this fiber.  The composition
\begin{equation}
\xymatrix{
\Omega(\FM{4}{I}\times \K)  \ar[d]^{\pi_*} &  &  \Omega(S^2\times S^2) 
\ar[ll]_(0.4){h_{13}^*\times h_{24}^*}\\
 \Omega\K   &  
}
\end{equation}
will thus produce a function on $\K$.

By Stokes' Theorem, the 
question of whether this function is a closed 0-form, i.e. a knot 
invariant, is now a question about the vanishing of the pushforward of $\alpha$ along 
the codimension one faces of $\FM{4}{I}$:
\begin{equation}\label{E:Stokes}
d\pi_{*}\alpha = \pi_{*}d\alpha - (\partial\pi)_{*}\alpha,
\end{equation}
where $(\partial\pi)_{*}\alpha$ means the sum of the pushforwards of 
$\alpha$ restricted to the codimension one faces.
However, since each $\omega_{ij}$ is a closed form, so is $\alpha$, and 
thus $\pi_{*}d\alpha =0$.  It follows that
\begin{equation}\label{E:SimplifiedStokes}
d\pi_{*}\alpha =  - (\partial\pi)_{*}\alpha.
\end{equation}

\begin{figure}[h]
\begin{center}
\input{Type2CD.pstex_t}
\caption{Diagram corresponding to integration along interior of 
$\FM{4}{I}$}\label{F:IntegrationDiagram}
\end{center}
\end{figure}
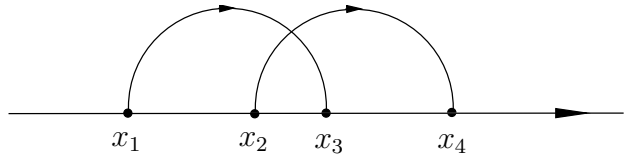

Since $(\partial\pi)_{*}\alpha$ is not zero, however, one is led to 
search for 
a ``correction term,'' namely another integral over a space with the 
same codimension one boundary as $\FM{4}{I}$.  A candidate presents 
itself immediately if we think of the first integral in terms of a 
chord diagram in Figure \ref{F:IntegrationDiagram}, where $x_{1}, \ldots, x_{4}$ represent the configuration points on the 
interval and the chords 
represent the way $h_{13}\times h_{24}$ pairs them off.

\begin{figure}[h]
\begin{center}
\input{Type2CDStratum.pstex_t}
\caption{Diagram corresponding to integration along the face 
$x_{2}=x_{3}$}\label{F:FaceDiagram}
\end{center}
\end{figure}
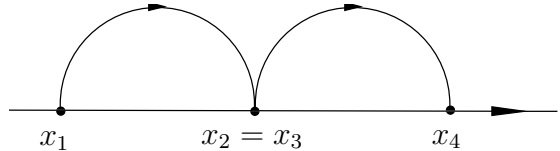

\noindent
The stratum $x_{2}=x_{3}$, for example, can then be pictured as 
another diagram as in Figure \ref{F:FaceDiagram}.

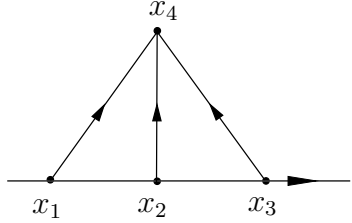
\begin{figure}[h]
\begin{center}
\input{Type2TD.pstex_t}
\caption{Simplest trivalent diagram}\label{F:SimpleTrivalent}
\end{center}
\end{figure}

Now consider Figure \ref{F:SimpleTrivalent}, the simplest
diagram with a trivalent vertex.
\noindent
If we let $x_{4}=x_{2}$ in that figure, we get exactly the same 
picture up to relabeling as the one depicting the stratum 
$x_{2}=x_{3}$ from above. 
The space one needs to study then turns out to be a compactified 
configuration space of four points in $\R^{3}$, three of which are 
restricted to 
lie on the knot.  Additionally, as the edges of the trivalent diagram suggests, 
there should 
be three maps to the sphere, $h_{12}$, $h_{13}$, and $h_{14}$, pulling 
back the volume forms.  
To construct the space itself, Bott and Taubes make the following general 
definition:

\begin{definition}\label{D:B-TPullback}  Define $F[k,s; \K, S^3]$ to be the pullback of
$$
\xymatrix{
                           &   \FM{k+s}{S^3}  \ar[d]   \\
\FM{k}{I}\times \K  \ar[r] &   \FM{k}{S^3}
}
$$
where the vertical map is the projection onto the first $k$ factors and 
the horizontal map is the evaluation of a knot in $\K$ on $k$ 
configuration points in $I$.
\end{definition}
The importance of these spaces is summarized in the following

\begin{prop}[\cite{BT}, Proposition A.3]\label{P:PullbackIsGood}
The space $F[k,s; \K, S^3]$ fibers 
over $\K$ and the fibers are compact manifolds with 
corners.  
\end{prop} 
A point in the fiber of 
$
F[k,s; \K, S^{3}] \longrightarrow  \K
$
may then be thought of as a $(k+3s)$-dimensional space of 
configurations of $k+s$ points in $S^{3}$, 
$k$ of which are restricted to lie on a given knot $K\in\K$.

Since one cannot consider vectors between 
configuration 
points in $S^{m}$, one instead replaces $S^{m}$ by 
$\R^{m}\cup\infty$, turning knots in $S^{m}$ into knots 
in $\R^{m}$.  This does not change any of the computations,
 except that one has to consider the \emph{strata at 
infinity}, or the boundary components of $F[k,s; \K, \R^{3}]$ 
determined by configuration points tending to infinity. 

Getting back to the trivalent diagram motivating \refD{B-TPullback}, the correction 
to the 
original integral $\pi_{*}\alpha$ turns out to be supplied by 
$F[3,1;\K,\R^{3}]$ as predicted.  More precisely, the edges in that 
diagram are thought of as
giving a prescription for pulling back the volume forms $\omega_{14}$, 
$\omega_{24}$, and $\omega_{34}$ from three spheres 
and then pushing forward to $\K$ the resulting 6-form  
$$
\alpha'=\alpha_{14}\alpha_{24}\alpha_{34}=h_{14}^{*}\omega_{14}\, 
h_{24}^{*}
\omega_{24}\, h_{34}^{*}\omega_{34}.
$$
Bott and Taubes show that the pushforwards of $\alpha$ and 
$\alpha'$ along all strata, including those at infinity, either vanish or 
cancel out.  We summarize in
the following theorem (also proved by Bar-Natan using 
different methods in \cite{BN3} and considered by Polyak and Viro in \cite{PV}):
\begin{thm}[\cite{BT}, Theorem 1.3]\label{T:Type2}  
The difference of 
pushforwards
$
\pi_*\alpha - \pi'_*\alpha'
$
is a knot invariant.
\end{thm}

Thurston generalizes in \cite{Th} as follows:
Suppose a 
labeled trivalent diagram $D$ is given and suppose it contains $k$ 
interval and $s$ free vertices.  The total number of its chords and 
edges is then $(k+3s)/2$.  These can be used as prescriptions for 
setting up as many maps 
$$
h_{ij}\colon F[k,s;\K,\R^{3}]\longrightarrow S^2.
$$
The volume form may be pulled back to 
$F[k,s;\K,\R^{3}]$ from each of the spheres, yielding a 
$(k+3s)$-form
\begin{equation}\label{E:Alpha}
\alpha = h_D^*\omega=\prod\limits_{\substack{\text{chords and} \\ \text{edges }ij}}
h_{ij}^*\omega_{ij}.
\end{equation}
Here $\omega$ is the product of the volume forms and $h_D$ the product of the $h_{ij}$.
Integration along the $(k+3s)$-dimensional 
fiber of $F[k,s;\K,\R^{3}]\to\K$
thus produces a function on $\K$, which we denote by $I(D,K)$.  Let $D_1$ be the one-chord diagram and let $k+3s=2n$.  We then have 
\begin{thm}[\cite{Alt, Th}]\label{T:Thurston}
 For a primitive weight system $W\in\W_{n}$, $n\geq 1$, the map
 $
 T(W)\colon\K\to \R
 $
 given by 
$$
T(W)(K)= \frac{1}{(2n)!}
\sum\limits_{D\in TD_n}
W(D)(I(D,K)-M_DI(D_1,K)),
$$
where $M_D$ is a real number which depends on $D$, is a knot invariant.
\end{thm}
This is also a finite type $n$ invariant, as was shown by Altschuler 
and Freidel \cite{Alt} (see \S\ref{S:Universal}).  The ``correction" term $M_DI(D_1,K)$ will be discussed in detail in \S\ref{S:AnomalousFaces}.

\begin{rems}  Recall that chord diagrams are also considered to be 
trivalent, so that the sum in \refT{Thurston} includes those as 
well.  Spaces $F[k,0;\K,\R^{3}]$ simply reduce to $\FM{k}{S^{1}}\times \K$
for those diagrams.

Next, fibers of $F[k,s;\K,\R^{3}]\to\K$ inherit their 
orientation from $\FM{k+s}{\R^{3}}$.  Trivalent 
diagrams are also oriented by the labeling of their vertices (edges 
and chords, as required in \S\ref{S:Trivalent}, always point from 
the vertex with the lower label, and this means that the 
vectors
 $h_{ij}$ must point in the corresponding directions).  Since a labeling 
 of a diagram determines the labeling of configuration points in 
 $F[k,s;\K,\R^{3}]$, changing the orientation of the diagram by permuting 
 vertex labels may change the orientation of $F[k,s;\K,\R^{3}]$.  The 
 corresponding integrals $I(D,K)$ and $I(D_1,K)$ will also differ in sign 
 for two labeled diagrams with different orientations.  But the weight 
 system $W$ also depends on the orientation, so that the signs will 
 cancel out.  Since the sum in the above theorem consists of $(2n)!$ identical 
 expressions, the normalizing factor of $1/(2n)!$ is introduced.


\end{rems}

We devote the next section to the proof of the above theorem.

\section{Vanishing of Integrals Along Faces}\label{S:Vanishing}

We prove \refT{Thurston} by checking that the integrals 
$(\partial\pi)_{*}\alpha$ on the codimension one strata of $F[k,s;\K,\R^3]$ 
either vanish or cancel out within the sum.  Throughout this section, different arguments are 
used for various types of faces, which are, using the terminology of 
Bott and Taubes, called
\begin{itemize}
\item \emph{principal}, if exactly two points 
degenerate,
\item \emph{hidden}, if more than two, but not all, points degenerate,
\item \emph{anomalous}, if all points degerate,
\item \emph{faces at infinity}, if one or more points approach 
infinity. 
\end{itemize}
  Using explicit
 coordinates on $F[k,s;\K,\R^{3}]$ we first show how $\alpha$ extends smoothly to such strata.  
\refT{Thurston} will follow as the combination of all the 
results in this section, at the end of 
\S\ref{S:AnomalousFaces}.  For the most part, the arguments used are given or suggested by Bott-Taubes and Thurston.


\subsection{Coordinates and pullbacks of forms on compactified configuration 
spaces}\label{S:Coordinates&Forms}

We first describe coordinates on manifolds $F[k,0;\K,\R^{3}]=F[k,I]\times\K$ and then
on the $F[k,s;\K,\R^{3}]$.   For the proof that we indeed get coordinates, see \cite{AS, FM}.
The evaluation map 
$$
F[k,I]\times\K \longrightarrow \FM{k}{\R^3}
$$
is given on the interior of $F[k,I]$ by 
$$
(x_{1}, \ldots, x_{k})\longmapsto p=(p_{1}, \ldots, p_{k})=
(K(x_{1}),\ldots, K(x_{k})).
$$
We now wish to extend this to the codimension one faces of $\FM{k}{I}$.  
So let $A$ be a subset of $\{1, \ldots, k\}$ containing $a$ 
elements, and let
 $q$ be a point in $\FM{k}{I}$ where all $x_{i}$ with $i\in A$ came 
 together at the same time.

Then $q$ is a configuration of the remaining $k-a+1$ 
points as well as a number of unit vectors recording the (one of the two 
possible, for each pair) directions of approach of the $a$ points.  
Applying $K$ 
to $q$ should thus yield a configuration of $k-a+1$ points in 
$\FM{k}{\R^{3}}$ 
as well as directions of approach of the colliding points
$K(x_{i})$.  We denote the stratum in 
$\FM{k}{\R^{3}}$ to which $q$ belongs by $\St{A}{\emptyset}$.

Now parametrize the neighborhood of the $a$ colliding points in $\FM{k}{I}$ by
\begin{equation}\label{E:StratumParameters-a0}
(x_{1}, u_{1}, u_{2}, \ldots, u_{a}, r\, ;\, x_{2}, \ldots, x_{k-a+1}),
\end{equation}
where
\begin{equation}\label{E:StratumConditions-a0}
x_{i}\in I \ \text{distinct}, \ \ \ u_{i}\in \R \ \text{distinct},\ \ \ r\geq 0, \ \ \ 
\sum_{i=1}^{a}u_{i}=0, \ \ \ \sum_{i=1}^{a}|u_{i}|^{2}=1.
\end{equation}
Then $q$ is given in the neighborhood of the stratum in question by
\begin{alignat*}{2}
q_{i} & =x_{1}+ru_{i}, \ \ \ &  i & \in\{1, \ldots, a\}, \\
q_{i} & =x_{i-a+1},    & i & \in\{a+1, \ldots, k\}.
\end{alignat*}
As $r$ approaches 0, we are in the limit left with the 
configuration $(x_{1}, x_{a+1}, \ldots, x_{k})$ as well as the direction 
vectors $u_{i}$.

\begin{rems} Notice that if two points on the interval, and later on 
a knot, collide, so will all points between them.

To simplify notation, we have chosen to label the 
colliding points by $1, \ldots, k$.
However, we do 
not require that configuration points are ordered in any standard way, 
so that all the indices we use should be considered up to 
permutation.
\end{rems}
This neighborhood in $\FM{k}{I}$ maps to a neighborhood 
of $\St{A}{\emptyset}$ in $\FM{k}{\R^{3}}$ consisting of points
\begin{alignat*}{2}
p_{i} & =K(x_{1}+ru_{i}), \ \ \ & i & \in\{1, \ldots, a\}, \\
p_{i} & =K(x_{i-a+1}), & i & \in\{a+1, \ldots, k\}.
\end{alignat*}
When $r=0$, the remaining configuration is
$$
(p_{1}, p_{a+1},\ldots, p_{k})=(K(x_{1}), K(x_{2}), \ldots, 
K(x_{i-a+1})),
$$
while the unit vectors recording the directions of the collision are taken to be 
$$
\lim_{r\to 0}\, \frac{K(x_{1}+ru_{j})-K(x_{1}+ru_{i})}
{|K(x_{1}+ru_{j})-K(x_{1}+ru_{i})|}, \ \ \ 
i,j\in\{1, \ldots,a\}, \ \ i<j.
$$
Upon expansion in Taylor series, this limit is 
easily seen to be
\begin{equation}\label{E:TangentialMap}
\tau =\frac{K'(x_{1})}{|K'(x_{1})|}.
\end{equation}
In analogy with (\ref{E:Tangentialmap}), we call $\tau$ the \emph{tangential map}.

Finally notice that the conditions (\ref{E:StratumConditions-a0}) also allow for the 
geometric interpretation of a point in $\St{A}{\emptyset}$ as an ``infinitesimal 
polygon modulo translation and scaling.''

\vskip 5pt
\noindent
To generalize to $F[k,s;\K,\R^{3}]$, take the coordinates on its interior to 
be
$$
(x_{1}, \ldots, x_{k}, x_{k+1}, \ldots, x_{k+s}),
$$
with
\begin{align*}
x_{i}\in I \text{ distinct}, & \ \ \ i\in\{1, \ldots, k\} \\
x_{i}\in \R^3 \text{ distinct}, & \ \ \ i\in\{k+1, \ldots, k+s\}.
\end{align*}
The points in the interior of $F[k,s; \K, \R^3]$ will then be given by
$$
(p_{1}, \ldots, p_{k+s})=(K(x_{1}), \ldots, K(x_{k}), 
x_{s+1}, \ldots, x_{k+s})
$$
if we additionally require 
$$
K(x_{i})\neq x_{j}, \ \ \ 1\leq i\leq k<j\leq k+s.
$$
This condition prevents against a point ``on the knot $K$'' coinciding with a 
point ``off the knot'' (in $\R^3$).  The terminology may be a bit misleading; keep in mind that a point ``off the knot'' 
is free to be anywhere in $\R^3$, including on the knot itself.

Now again assume for simplicity that the points on the knot are indexed 
by $1, \ldots, k$, and those off the knot by $k+1, \ldots, k+s$.  
Let $A$ and $B$ be subsets of $\{1, \ldots, k\}$ and $\{k+1, 
\ldots, k+s\}$ with cardinalities $a$ and $b$.  Denote by $\St{A}{B}$ 
the stratum of $F[k,s; \K, \R^3]$ given by the 
coming together of the $a+b>1$ points indexed by the elements of $A$ 
and $B$.  Assume for now that $A$ is nonempty, in which case the limit point is necessarily on the knot.  A remark on the case $A=\emptyset$ follows \eqref{E:StratumDirections-ab}. 

For each $A$ and $B$, we thus get a stratum $\St{A}{B}$.  However, 
some of the strata are empty because if two points on the knot 
collide, so will all points between them.  Thus if $A$ contains 
indices $i$ and $j$, but not all the indices of points between 
$p_{i}$ and $p_{j}$, $\St{A}{B}$ is empty.

To give coordinates on a neighborhood of $\St{A}{B}$, assume as 
before that $K$ is given and introduce 
parameters
\begin{equation}\label{E:StratumParameters-ab}
(x_{1}, u_{1}, \ldots, u_{a+b}, r\, ;\, x_{2}, \ldots, x_{k+s-a-b+1})
\end{equation}
which satisfy the following conditions:
\begin{align}\label{E:StratumConditions-ab}
\text{(1)} & \ \ \ x_{i}\in I  \text{ distinct}, \ \ 
1\leq i\leq k-a+1, \notag \\
\text{(2)} & \ \ \ x_{i}\in \R^3 \text{ distinct}, \ \ 
k-a+2\leq i\leq k+s-a-b+1,   \notag \\
\text{(3)} & \ \ \   u_{i}\in \R \text{ distinct}, \ \ 
1\leq i\leq a, \notag \\ 
\text{(4)} & \ \ \   u_{i}\in \R^3 \text{ distinct}, \ \ 
a+1\leq i \leq a+b, \notag \\
\text{(5)} & \ \ \   r\geq 0, \\
\text{(6)} & \ \ \  K(x_{i})\neq x_{j}, \ \ 
i\leq k-a+1, \ \ j<k-a+1, \notag \\
\text{(7)} & \ \ \  K'(x_{1})u_{i}\neq u_{j}, \ \ 
i\leq a, \ \ j>a \notag \\
\text{(8)} & \ \ \   
\sum_{i=1}^{a}|K'(x_{1})|^{2}u_{i}^{2}+\sum_{i=a+1}^{a+b}|u_{i}|^{2}=1
\notag \\
\text{(9)} & \ \ \  
\sum_{i=1}^{a}u_{i}+\sum_{i=a+1}^{a+b}
\frac{\langle K(x_{1}), u_{i}\rangle}{|K(x_{1})|^{2}}=0. 
\notag
\end{align}
Conditions (\ref{E:StratumConditions-ab}.8) and 
(\ref{E:StratumConditions-ab}.9) ensure that the limiting directions 
between colliding points are indeed vectors in the unit sphere bundle 
of the normal bundle, and can again be thought of as scaling 
and translation of vectors in $\R^3$.  Condition 
 (\ref{E:StratumConditions-ab}.6) prevents against a point on the 
 knot colliding with a point off the knot before the rest of the 
 points join them, so that, along with 
 (\ref{E:StratumConditions-ab}.1) and 
 (\ref{E:StratumConditions-ab}.2), the stratum $\St{A}{B}$ is 
 described by the $a+b$ points coming together exactly at the same time.

Configuration points near $\St{A}{B}$ ($r>0$) are then given by
\begin{alignat}{3}\label{E:StratumPoints-ab}
p_{i} & =K(x_{1}+ru_{i}), \ \ \ & i & \in\{1, \ldots, a\} &
\ \  & \text{(on knot, colliding)},  \notag \\
p_{i} & =K(x_{i-a+1}), & i & \in\{a+1, \ldots, k\} &
\ \  & \text{(on knot, not colliding)},  \\
p_{i} & =K(x_{1})+ru_{i-k+a}, \ \ \ & i & \in\{k+1, \ldots, 
k+b\} &
\ \  & \text{(off knot, colliding)},  \notag \\
p_{i} & =x_{i-a-b+1}, & i & \in\{k+b+1, \ldots, k+s\} &
\ \  & \text{(off knot, not colliding)}. \notag
\end{alignat}

Again remember that the indexing of points depends on the labeling of 
the configuration points as well as the sets $A$ and $B$.

When $r=0$, what is left is a configuration of $k+s-a-b+1$ points with 
the colliding points becoming $K(x_{1})$.  The limiting directions 
are

\begin{alignat}{2}\label{E:StratumDirections-ab} 
\tau =\frac{K'(x_{1})}{|K'(x_{1})|}, &   \ \ \ & & \text{$p_{i}$, $p_{j}$ on the knot},  
\notag \\ 
\frac{u_{j}-K'(x_1)u_{i}}{|u_{j}-K'(x_1)u_{i}|}, &  \ \ \  & & \text{$p_{i}$ on the knot, $p_{j}$ off the 
knot},  \\
\frac{u_{j}-u_{i}}{|u_{j}-u_{i}|}, &  \ \ \ & & \text{$p_{i}$, $p_{j}$ off the knot}.  \notag
\end{alignat}
Constraint 
(\ref{E:StratumConditions-ab}.7) prevents against the second vector 
being 0.

It should be clear how the above has to be modified in the case $A=\emptyset$.  The colliding points are all off the knot know, so in particular the conditions (\ref{E:StratumParameters-ab}) simplify.  Some are vacuous, and condition (9) has to be changed to reflect the fact that $x_1$ should be a point in $\R^3$ off the knot ``in the middle" of the colliding points. Same is true for the third equation in \eqref{E:StratumPoints-ab}, while the first is not needed at all.

\vskip 5pt
\noindent
We are interested in how $\alpha=h_D^*\omega$ from \eqref{E:Alpha} 
 restricts to 
each $\St{A}{B}$.  For this, it suffices to  know how the maps 
$h_{ij}$ behave on $\St{A}{B}$.  The answer is now simple:  If 
$h_{ij}$ pairs two points, one 
or both of which are not among the colliding points, then it restricts 
to $\St{A}{B}$ as the normalized difference of two points in 
(\ref{E:StratumPoints-ab}), with $r=0$ (which two depends 
on whether the points are on or off the knot).  If both points are 
among the colliding ones, $h_{ij}$ extends to $\St{A}{B}$ as one of the  
maps in (\ref{E:StratumDirections-ab}).  Since these extensions are 
smooth, the pullback $\alpha$ also restricts smoothly to those 
codimension one strata given by collisions of points.  

The remaining 
codimension one strata are given by one or more points going to 
infinity.  We parametrize such strata the same way, except $r$ now 
tends to infinity rather than 0.  The limiting map between two 
points off the knot, $p_{i}$ and $p_{j}$, going to 
infinity is then 
\begin{equation}\label{E:InfinityDirection1}
h_{ij}=\lim_{r\to\infty}
\frac{(x_{1}+ru_{j})-(x_{1}+ru_{i})}{|(x_{1}+ru_{j})-(x_{1}+ru_{i})|}=
\frac{u_{j}-u_{i}}{|u_{j}-u_{i}|},
\end{equation}
which is exactly the same as the last map in 
(\ref{E:StratumDirections-ab}), so nothing important changes in this case either.
The 
map between a point $p_{j}$ going to infinity and a point $p_{i}$ 
which does not, also extends smoothly to the stratum at infinity as
\begin{equation}\label{E:InfinityDirection2}
h_{ij}=\lim_{r\to\infty}
\frac{(x+ru_{j})-p_{i}}{|(x+ru_{j})-p_{i}|}=\frac{u_{j}}{|u_{j}|}.
\end{equation}

In the following sections, we closely examine these maps to show that, 
for each $A$ and $B$, the 
pushforward of $\alpha$ along $\St{A}{B}$ vanishes, cancels with 
others, or has to be compensated for by another integral.

\subsection{Faces determined by disconnected sets of vertices}\label{S:DisconnectedStratum}

Given a labeled trivalent diagram with $k$ interval and $s$ free 
vertices, we now think of its vertex set as determing a space $F[k,s; \K, \R^3]$, 
and of its chords and edges as determining the map $h_{D}$.  The 
strata $\St{A}{B}$ can also be thought of as being determined by 
subsets $A$ and $B$ of the vertex set since a labeling of the diagram 
determines a labeling of the configuration points in $F[k,s; \K, \R^3]$.  
Elements of $A$ are interval vertices, while the vertices listed in 
$B$ are free.

Suppose then a stratum is represented by a subset $A\cup B$ of a 
diagram's vetex set, and also suppose the chords and edges 
of the diagram are such that $A\cup B$ can be broken up into at least 
two smaller 
subsets such that no chord or edge connects a vertex of one subset to 
a vertex of another.  We will say that such a set of vertices $A\cup B$
is \emph{disconnected}, and 
denote the stratum corresponding to $A\cup B$
by $\Sd$.

\begin{prop}\label{P:disconnected}  Unless $a=2$ and $b=0$, the 
pushforward 
of $\alpha$ to $\K$ along $\Sd$ vanishes.
\end{prop}
The reason this statement is true is essentially that the two connected components can be translated independently.  However, now write down the proof in some detail since we will refer to it often throughout.  We first need the following lemma whose 
proof is left to the reader:
\begin{lemma}\label{L:LowerDimension}
Suppose $X$ and $Y$ are spaces which fiber over $Z$ with compact 
fibers of dimensions $m$ and $n$, and suppose $f$ is a fiber-preserving map 
from $X$ to $Y$.  Let $\beta$ be a $p$-form on $Y$ and let 
$\gamma=f^{*}\beta$.  If 
$m>n$, then the pushforward of $\gamma$ from $X$ to 
$Z$ vanishes.
\end{lemma}

\begin{proof}[Proof of \refP{disconnected}]  Given a trivalent 
diagram $D$ and a disconnected subset of vertices of $D$ determining a
stratum $\Sd$, we will construct a space $\mathcal{S}'$ of dimension 
strictly 
lower then the dimension of $\Sd$ such that the map $h_{D}$ 
factors through $\mathcal{S}'$:
\begin{equation}
\xymatrix{
\Sd \ar[rr]^{h_{D}} \ar[dd]^{\partial\pi} \ar[dr]^{f}
 &                       &   (S^{2})^{|e|}\\
 &  \mathcal{S}'\ar[dl]^{\partial\pi'} \ar[ur]^{h_{D}'} &   \\
\K
}
\end{equation}

Suppose first that $\Sd$ corresponds to a set  of vertices $A\cup B$ with 
 exactly two disconnected
subsets.  Let there be $a_{1}$ and $a_{2}$ interval vertices in 
each of those subsets, respectively, with $a_{1}+a_{2}=a=|A|$.  Similarly, 
let
$b_{1}+b_{2}=b=|B|$, 
where $b_{1}$ and $b_{2}$ are the numbers of free vertices in each 
subset.  Assume without 
loss of generality that labeling of $D$ is such that the $a$ 
vertices are labeled by $1, \ldots, a$ and the rest by $a+1, \ldots, 
a+b$.

The parameters 
describing a 
neighborhood of the colliding points in $\Sd$ are 
$$(x, u_{1}, u_{2}, \ldots, u_{a+b}, r),$$
where $u_{1}, ..., u_a$ are in $I$, and $u_{a+1}, ..., u_{a+b}$ are 
in $\R^3$.  For now, we omit the 
parameters for configuration points in $\Sd$ which do not 
come together.  We will get back to them later. 

The $u_{i}$ also must 
safisfy conditions (\ref{E:StratumConditions-ab}) (with $x_{1}$ now 
replaced by $x$ to simplify notation).  Point $x$ is in $I$ 
unless $a=0$, 
in which case $x$ is a vector in $\R^3$.

Given $K\in\K$, the colliding points are
\begin{align}
p_{i}=K(x+ru_{i}), \ \ \  & i\in\{1, \ldots, a\}\notag \\
p_{i}=K(x)+ru_{i}, \ \ \  & i\in\{a+1, \ldots, a+b\}, \notag 
\end{align}
as before.  Recall also that the maps $h_{ij}$ are given on $\Sd$ by
\eqref{E:StratumDirections-ab}.

Now let $F[k,s; \K, \R^3]'$ be a space differing from $F[k,s; \K, \R^3]$ in that the 
blowups are not performed along all the diagonals involving points 
 at least 
one of which corresponds to a vertex in one subset of $A\cup B$
and at least one 
of which corresponds to the other.  Instead, the diagonals are simply put back in so that now not all directions of approach of points are 
recorded.  This is still a compact manifold with 
corners since Fulton-MacPherson compactification can be constructed by a sequence of blowups along the diagonals and one obtains a manifold with corners at each stage \cite{FM, AS}.  Thus $F[k,s; \K, \R^3]'$ is in fact a submanifold with corners of $F[k,s; \K, \R^3]$.

A part of the boundary of $F[k,s; \K, \R^3]'$ 
consists of the same $a+b$ points as in $\Sd$ coming together.  
Denote this face by $\mathcal{S}'$ and
consider its neighborhood described by two sets of 
independent parameters
\begin{equation}\label{E:DisconnectedParameters}
(x_{1}, v_{1}, \ldots, v_{a_{1}+b_{1}}, r_{1}\, ;\, 
x_{2}, , v_{a_{1}+b_{1}+1}\ldots, v_{a+b}, r_{2}).
\end{equation}
If $a_1>0$, then $x_1$ is in $I$ and it is in $\R^3$ if $a_1=0$.  The 
same holds for $a_2$ and $x_2$.  If both $x_{1}$ and $x_{2}$ are in 
$I$, we may assume $x_{1}<x_{2}$.
Each set of $v_i$ must satisfy the conditions 
(\ref{E:StratumConditions-ab}).  Thus in particular, for the first set of parameters, we 
have
\begin{equation}\label{E:DisconnectedConstraint1}
\sum_{i=1}^{a_{1}}|K'(x_{1})|^{2}v_{i}^{2}+
\sum_{i=a_{1}+1}^{a_{1}+b_{1}}|v_{i}|^{2}=1,
\end{equation}
\begin{equation}\label{E:DisconnectedConstraint2}
\sum_{i=1}^{a_{1}}v_{i}+
\sum_{i=a_{1}+1}^{a_{1}+b_{1}}\frac{\langle K'(x_{1}),v_{i} 
\rangle}{|K'(x_{1})|^{2}}=0.
\end{equation}
As usual, these parameters describe $a_1+b_1$ points 
\begin{align}
p_i'=h_t(x_1+r_1v_i), \ \ \ & i\in\{1, ..., a_1\},\notag \\
p_i'=h_t(x_1)+r_1v_i, \ \ \ & i\in\{a_1+1, ..., a_1+b_1\}.\notag
\end{align}
The restrictions of the maps $h_{ij}'$ to $\mathcal{S}'$ are identical 
to those for $\Sd$
up to renaming of the parameters..

We now proceed to construct $f$ for this subset of $A\cup B$.
If $x$ is a vector in $\R^3$, then $x_{1}$ will be as well, so set 
$x_{1}=x$.  The same happens if $x$ is in $I$ and the subset
 has an interval vertex.  However, it may 
happen that $x$ is in $I$ (this means that $A\cup B$ contains at least 
one interval vertex), while $x_{1}$ is in 
$\R^3$ (this means that the subset of $A\cup B$
has no interval vertices).  In this case, 
we set $x_{1}=K(x)$.  We 
also set $r_{1}=r$.
 
Remember that we ultimately want
\begin{equation}\label{E:CompositionPreserved}
h_{ij}'\circ f=h_{ij}
\end{equation}
when $r_{1}=r=0$.  Note that $h_{ij}$ can only be the tangential map  
(first map in \eqref{E:StratumDirections-ab}) when $A\cup B$ contains two
interval vertices which are connected by a chord.  But these two 
vertices then form a subset of $A\cup B$ which is 
disconnected from the rest of $A\cup B$.  The first set of 
parameters in \eqref{E:DisconnectedParameters} can then be taken as
$(x_{1}, r_{1})$, $x_{1},r_{1}\in I.$
The only requirement $f$ now has to satisfy is thus 
$$\frac{K'(x_1)}{|K'(x_1)|}\circ 
f=\frac{K'(x)}{|K'(x)|},$$
and this is true since $f$ sends $x$ to $x_{1}$.  An important 
consequence of this observation is \refC{anomalous-cd}.


We can thus assume the first subset of $A\cup B$ does not consist of 
two vertices with a chord between them, and we turn our attention to 
the remaining two maps in (\ref{E:StratumDirections-ab}).

Suppose $f$ gives $v_{i}=z_{i}$ for some $z_{i}$, 
$1\leq i\leq a_{1}+b_{1}$ ($a_{1}$ of these are numbers and $b_{1}$ 
are vectors).    In order for \eqref{E:CompositionPreserved} to hold 
for the last two maps in \eqref{E:StratumDirections-ab},
 the $z_{i}$ should satisfy
\begin{align}
z_{j}-K'(x)z_{i}=u_{j}-K'(x)u_{i}, \ \ \ 
&  i\in\{1, \ldots, a_{1}\}, \ j\in\{a_{1}+1, \ldots, 
a_{1}+b_{1}\} \label{E:zCondition1} \\
z_{j}-z_{i}=u_{j}-u_{i}, \ \ \ 
&  i,j\in\{a_{1}+1, \ldots, a_{1}+b_{1}\}. \label{E:zCondition2}
\end{align}
Since the vertices in the subset of $A\cup B$ we are considering are 
connected by edges, it is easily seen that there will now be 
exactly one fewer 
independent equations then functions $z_i$.  
But these must also satisfy the constraints 
(\ref{E:DisconnectedConstraint1}) and (\ref{E:DisconnectedConstraint2}).  The 
second constraint can be added to the system 
\eqref{E:zCondition1}--\eqref{E:zCondition2} which will now 
produce a unique solution.  The $z_i$ will be various
combinations of vectors $u_{j}-u_{i}$, $u_{j}-K'(x)u_{i}$,
$K'(x)$, and their magnitudes.  They already satisfy the 
constraint (\ref{E:DisconnectedConstraint2}), and to make sure they also 
satisfy (\ref{E:DisconnectedConstraint1}), we may simply 
divide each $z_{i}$ by
$$\left(\sum_{i=1}^{a_{1}}|K'(x_{1})|^{2}z_{i}^{2}+
\sum_{i=a_{1}+1}^{a_{1}+b_{1}}|z_{i}|^{2}\right)^{\frac{1}{2}}.$$
The above expressions are never zero since
\begin{align}
u_{j}\neq K'(x)u_{i}, \ \ \ & i\in\{1, \ldots, a_{1}\}, \ 
j\in\{a_{1}+1,\ldots, a_{1}+b_{1}\},\notag \\
u_{j}\neq u_{i}, \ \ \ & i,j\in\{a_{1}+1,\ldots, a_{1}+b_{1}\}.\notag
\end{align}
Modifying the $z_{i}$ by these factors does not affect the compositions 
$h_{ij}'\circ f$ either.  The factors cancel in $h_{ij}'$ since they 
are 
positive, real-valued functions.  The compositions are thus still 
$h_{ij}$ 
as desired.

This procedure can be repeated to construct $f$ for the other set of 
parameters 
in $\mathcal{S}'$, where $f$ should set $x_{2}=x$ or 
$x_{2}=K(x)$, and 
$r_{2}=r$.

For simplicity, the coordinates of the points in $\Sd$ that do not 
come together have 
been 
omitted from the previous discussion. But 
the coordinates for those points are the same on $\Sd$ and 
$\mathcal{S}'$, so $f$ is the identity there, and it immediately follows that
$h_{ij}'\circ f=h_{ij}$ for all the possible cases of one or both points outside of $\Sd$.

We have thus constructed a space $\mathcal{S}'$ through which $h_{D}$ 
factors.  
The parameters in \eqref{E:DisconnectedParameters}, 
describing a neighborhood of $\mathcal{S}'$, determine, along with 
the constraints for 
each set of them, precisely as many dimensions as those 
describing a neighborhood of $\Sd$ in $F[k,s; \K, \R^3]$.  But $f$ then 
restricts
 $x_{1}$, $x_{2}$, $r_{1}$, and $r_{2}$ and  
it follows that the fiber dimension of $\mathcal{S}'$ is smaller than 
that of 
$\Sd$.

\refL{LowerDimension} thus completes the argument for this case of exactly 
two subsets of $A\cup B$ which have no chords or edges connecting them.  
The case of 
more 
subsets can be reduced to this situation by parametrizing the 
configuration points corresponding to all but 
one 
subset with one set of parameters.  We are thus 
viewing some number of subsets as two,
with one of them itself containing further subsets which are 
disconnected from each other.  The only change 
in the argument is that there will now be fewer maps $h_{ij}$ and 
$h_{ij}'$ 
for this subset
than before, and hence fewer conditions \eqref{E:zCondition1} and 
\eqref{E:zCondition2}. The $z_{i}$ 
will therefore not be unique.  The number of maps 
$h_{ij}$ is also the reason why this procedure does not work for a 
stratum determined by a single subset $A\cup B$ of the vertices of 
$D$ which cannot be broken down into more subsets not joined by 
edges:  If we try to reparametrize such a stratum by 
two sets of parameters, there will be 
more 
independent equations coming from the $h_{ij}$ then variables $z_{i}$.
\end{proof}

It is easy to see
why this proposition fails in the case 
 of a principal face 
corresponding to a subset of the vertex set of $D$ consisting of 
exactly two interval
vertices which are not connected by a chord.  The reason 
is essentially that, when two points in 
$\R^{2}$ 
come together at the blownup diagonal, the dimension of the 
resulting space is one, as it would have been had the blowing up not 
been performed, but rather the diagonal in $\R^{2}$ was put back in. 
This case, however, is 
taken care of in \S\ref{S:PrincipalFaces}.
\vskip 5pt
\noindent
An easy consequence of \refP{disconnected} is as follows:  Consider 
a diagram $D$ with at least one chord.  Then the vertex set of $D$ 
is necessarily disconnected, with the two vertices joined by a chord 
forming a subset of the vertex set which is not connected to other 
vertices.  Let 
$\mathcal{S}_D$ denote the anomalous face of $F[k,s; \K, \R^3]$, or 
the part of the boundary of $F[k,s; \K, \R^3]$ with all 
$k+s$ points coming together.  We then immediately have
\begin{cor}\label{C:anomalous-cd}   
If $k+s>2$ and a diagram $D$ contains a chord, then the pushforward 
of $\alpha=(h_D)^*\omega$ to $\K$ along $\mathcal{S}_D$ vanishes.
\end{cor}
In particular, if $D$ is a chord diagram (no free vertices), then the 
pushforward along 
$\mathcal{S}_{D}$ must be 0.  (In fact, the same reasoning says that 
the pushforward along any hidden face for a chord diagram is zero, 
but this case is covered in the next section.)  \refC{anomalous-cd} will be used in \refP{anomalous}.

\subsection{Hidden faces}\label{S:HiddenFaces}

In this section, we prove

\begin{prop}\label{P:hidden}  Suppose $\St{A}{B}$ is as before determined by the 
subset $A\cup B$ of the vertex set of a prime trivalent diagram $D$ with 
$2n$ vertices.  Suppose also $2<|A\cup B|<2n$.  Then the 
pushforward 
of $\alpha$ to $\K$ along $\St{A}{B}$ vanishes.
\end{prop}
We first need two lemmas which are essentially due to
  Kontsevich \cite{Kont3}.

\begin{lemma}\label{L:case2hidden}
Suppose $A\cup B$ contains vertices $i_{1}$, \ldots, 
$i_{4}$ with edges between them as in Figure \ref{F:Lemma1}.
Each of the vertices $i_{3}$ and $i_{4}$ is either on the interval 
or free, and the two edges emanating upward from $i_{1}$ 
end in vertices not in $A\cup B$.  Then the 
pushforward 
of $\alpha$ to $\K$ along $\St{A}{B}$ vanishes.
\end{lemma}
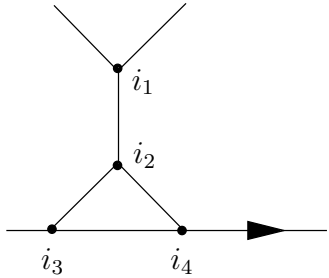
\begin{figure}[h]
\begin{center}
\input{Orientation1.pstex_t}
\caption{The case of Lemma \ref{L:case2hidden}}\label{F:Lemma1}
\end{center}
\end{figure}
\noindent

\begin{proof}
We will prove the statement by exhibiting an 
automorphism $\phi$
of $\St{A}{B}$ which preserves its orientation but takes $\alpha$ 
to $-\alpha$.  The three maps corresponding to the pictured edges will be permuted, but they will all also be negated.  Since the two pushforwards have to be the same,
it will follow that $(\partial\pi)_{*}\alpha$ equals 
its negative and hence must be 0.  

Recall that a neighborhood of $\St{A}{B}$ is parametrized by
$$(x_{1}, u_{1}, \ldots, u_{a}, u_{a+1}, \ldots, u_{a+b}, r\, ;\,  
x_{2}, 
\ldots, x_{k+s-a-b+1}).$$
Assume first that both vertices $i_{3}$ and $i_{4}$ are free and
consider a map $\phi$ from $\St{A}{B}$ to itself which leaves $r$, 
$x_{i}$ alone and sends $u_{i}$ 
to $w_{i}$ given by
\begin{alignat}{2}
w_{i_{1}}&= 
g_{1}(u_{i_{3}}+u_{i_{4}}-u_{i_{1}}+g_{2}K'(x_{1})), \ \ \ & & 
\label{E:wi1} \\
w_{i_{2}}&= 
g_{1}(u_{i_{3}}+u_{i_{4}}-u_{i_{2}}+g_{2}K'(x_{1})), & &\\
w_{i}&= 
g_{1}(u_{i}+g_{2}), \ \ \ \ \ & i &\in\{1, \ldots, a\},\\
w_{i}&= 
g_{1}(u_{i}+g_{2}K'(x_{1})), \ \ \ \ \ & i &\in\{a+1, \ldots, 
a+b\}\setminus \{i_{i}, i_{2}\},
\end{alignat}
and where
\begin{align}
g_{1}=&\bigg( \sum_{i=1}^{a}|K'(x_{1})|^{2}(u_{i}+g_{2})^{2}
              +\sum_{\stackrel{i=a+1}{i\neq i_{1},i_{2}}}^{a+b}
              |u_{i}+g_{2}|^{2} \\
      &  \ \  +|u_{i_{3}}+u_{i_{4}}-u_{i_{1}}+g_{2}K'(x_{1})|^{2}
              +|u_{i_{3}}+u_{i_{4}}-u_{i_{2}}+g_{2}K'(x_{1})|^{2}
       \bigg)^{-\frac{1}{2}}, \notag \\
g_{2}=&  -\frac{2\, \langle K'(x_{1}), u_{i_{3}}+u_{i_{4}}-u_{i_{1}}-u_{i_{2}}  \rangle}
          {(a+b)\, |K'(x_{1})|^{2}}.    \label{E:g2}     
\end{align}

To check that $\phi$ is indeed an automorphism, we need to show that 
it preserves conditions \eqref{E:StratumConditions-ab}.  All $w_{i}$, $i\in\{1, \ldots, a\}$, are 
distinct numbers as the corresponding $u_{i}$ are.  The same is true for 
vectors $w_{i}$, $i\in\{a+1, \ldots, a+b\}\setminus \{i_{i}, i_{2}\}$.  
Also, $w_{i_{1}}$ is different from $w_{i_{2}}$, but one or both of them 
might equal some other vectors.  

So suppose, for instance, that 
$w_{i_{1}}$ equals $w_{i}$ for some $i$.  In this case, we simply go 
back to the construction of $F[k,s; \K, \R^3]$ and modify it so that the blowup along the 
diagonal $p_{i_{i}}=p_{i}$ is not performed.  The effect on the 
parametrization of $\St{A}{B}$ is that the
vectors $u_{i_{1}}$ and $u_{i}$ are no longer required to be distinct.
Now the two points are allowed to come 
together before the others in $\St{A}{B}$ but, since there was no blowup 
along 
their diagonal, $\St{A}{B}$ still has codimension one 
in $F[k,s; \K, \R^3]$.  The important observation is that
this has no effect on 
the maps $h_{ij}$---there is no edge connecting $i_{i}$ 
and $i$ in $D$ and hence no map relating the corresponding points, 
$p_{i_{i}}$ and $p_{i}$.  (This kind of a modification in the construction of 
$F[k,s; \K, \R^3]$ was also required in the proof of \refP{disconnected}.)

By omitting some of the blowups, we can similarly preserve the condition
$$K'(x_{1})u_{i}\neq u_{j},\ \ \ i\in\{1, \ldots, a\}, \ j\in\{a+1, 
\ldots, a+b\}.$$
Now we do not impose this condition for $j=i_{1}$ and $j=i_{2}$ while, for the 
remaining $j$, it is still satisfied after $\phi$ is applied.

The condition
$$K(x_{i})\neq x_{j}, \ \ \ i\in\{1, \ldots, k-a+1\}, \ j\in\{k-a+2, 
\ldots, k+s-a-b+1
\},$$
trivially remains satisfied since $\phi$ is the identity on those 
parameters.

As for the remaining two conditions, 
constraint (\ref{E:StratumConditions-ab}.8) maps under $\phi$ to
\begin{align}
& \sum_{i=1}^{a}|K'(x_{1})|^{2}( g_{1}(u_{i}+g_{2})) ^{2}
              +\sum_{\stackrel{i=a+1}{i\neq i_{1},i_{2}}}^{a+b}
              |g_{1}(u_{i}+g_{2})|^{2} \notag \\
             & + |g_{1}(u_{i_{3}}+u_{i_{4}}-u_{i_{1}}+g_{2}K'(x_{1}))|^{2}
               +|g_{1}(u_{i_{3}}+u_{i_{4}}-u_{i_{2}}+g_{2}K'(x_{1}))|^{2}
               =\ 1 \notag 
\end{align}
Constraint (\ref{E:StratumConditions-ab}.9) becomes
\begin{align}
g_{1}\ \bigg( & \sum_{i=1}^{a}(u_{i}+g_{2})
                +\sum_{\stackrel{i=a+1}{i\neq i_{1},i_{2}}}^{a+b}
                \frac{\langle K'(x_{1}), u_{i}+g_{2}K'(x_{1})  \rangle}
                     {|K'(x_{1})|^{2}}  \notag \\
              & +\frac{\langle K'(x_{1}), 
                         u_{i_{3}}+u_{i_{4}}+u_{i_{1}}+g_{2}K'(x_{1})  \rangle}
                     {|K'(x_{1})|^{2}}
                +\frac{\langle K'(x_{1}), 
                         u_{i_{3}}+u_{i_{4}}+u_{i_{2}}+g_{2}K'(x_{1})  \rangle}
                     {|K'(x_{1})|^{2}}
         \bigg) 
         =\ 0 \notag 
\end{align}
Thus $g_{1}$ and $g_{2}$ are simply correction functions which ensure 
that $\phi$ preserves the last two requirements of 
(\ref{E:StratumConditions-ab}).

The automorphism $\phi$ only affects three of the maps $h_{ij}$,  
where the compositions $\phi\circ h_{ij}$ on $\St{A}{B}$ give
\begin{equation*}
\phi\circ h_{i_{1}i_{2}} 
                =-h_{i_{1}i_{2}}, \ \ 
\phi\circ h_{i_{2}i_{3}}  
                =-h_{i_{2}i_{4}}, \ \ 
\phi\circ h_{i_{2}i_{4}}  
                =-h_{i_{2}i_{3}}.
\end{equation*}
Therefore $\omega$ pulls 
back to $\St{A}{B}$ as $-\alpha$ (switching the two maps, 
$h_{i_{2}i_{3}}$ and $h_{i_{2}i_{4}}$, preserves $\alpha$ since this 
has the effect of switching two $(m-1)$-forms $\omega_{ij}$ with $m$ 
odd).

Also, $\phi$ does not change the orientation of $\St{A}{B}$.  
The easiest way to see this is to think of the parameters 
$u_{i}$, and consequently $w_{i}$, as numbers and vectors modulo translation and 
scaling as explained in section \S\ref{S:Coordinates&Forms}.  Functions $g_{1}$ 
and $g_{2}$ can then be dropped from the definition of $\phi$, whose 
Jacobian is then easily seen to have positive determinant.

A slight modification is required in the case that one or both of 
$i_{3}$ and $i_{4}$ are vertices on the interval.  In case 
$i_{3}$ is on the interval, $u_{i_{3}}$ has to be replaced by 
$K'(x)u_{i_{3}}$ in (\ref{E:wi1})--(\ref{E:g2}).  Same for $i_{4}$.  The rest of 
the argument is unchanged, and this completes the proof.
\end{proof}
\begin{lemma}\label{L:case1hidden} Suppose $A\cup B$ contains vertices 
$i_{1}$, $i_{2}$, and $i_{3}$ with edges between them as in Figure 
\ref{F:Lemma2}.
Each of the vertices $i_{2}$ and $i_{3}$ is either on the interval or free, 
and the third edge emanating from $i_{1}$ 
ends in a vertex not in $A\cup B$.  Then the 
pushforward 
of $\alpha$ to $\K$ along $\St{A}{B}$ vanishes.
\end{lemma}
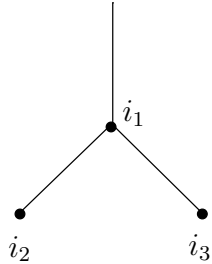
\begin{figure}[h]
\begin{center}
\input{Orientation2.pstex_t}
\caption{The case of Lemma \ref{L:case1hidden}}\label{F:Lemma2}
\end{center}
\end{figure}
\noindent

\begin{proof}  The argument is essentially the same as in the previous 
lemma.  The automorphism $\phi$ is 
given by
\begin{align}
u_{i_{1}}&\longmapsto 
g_{1}(u_{i_{2}}+u_{i_{3}}-u_{i_{1}}+g_{2}K'(x)), \ \ \ 
\notag\\
u_{i}&\longmapsto g_{1}(u_{i}+g_{2}), \ \ \ & i&\in\{1, \ldots, a\},\notag\\
u_{i}&\longmapsto g_{1}(u_{i}+g_{2}K'(x)), \ \ \ & i&\in\{a+1, \ldots, 
a+b\}\setminus i_{i}.\notag
\end{align}
As before, we may need to 
multiply one or both of $u_{i_{2}}$ and $u_{i_{3}}$ by $K'(x)$ 
depending on whether vertices 
$i_{2}$ and $i_{3}$ are on the interval or free.  Correction functions $g_{1}$ 
and $g_{2}$ are again defined so that the parameter
constraints are preserved.  Composing with the maps $h_{ij}$ gives
\begin{equation*}
\phi\circ h_{i_{1}i_{2}}=-h_{i_{1}i_{3}},\ \ 
\phi\circ h_{i_{1}i_{3}}=-h_{i_{1}i_{2}},
\end{equation*}
and $\alpha$ thus remains unchanged.  However, this 
automorphism reverses the orientation of the stratum, so that the 
pushforward of $\alpha$ along $\St{A}{B}$ again must be zero.  
\end{proof}

\begin{proof}[Proof of \refP{hidden}]
Recall that every trivalent diagram $D$ we are considering is prime, 
namely not a connected sum of two or more diagrams.  
It follows that all hidden faces $\St{A}{B}$ of 
$F[k,s; \K, \R^3]$ correspond to subsets of vertices of $D$ with at least one chord or edge 
connecting a vertex in $A\cup B$ with a vertex not in $A\cup B$. 
Let $v$ be such a vertex in $A\cup B$ and let $e_{v}$ be the number of
 chords 
or edges connecting it to vertices in the complement of $A\cup B$.  
 Then consider 
the following cases: 
\vskip 5pt
\noindent
$v \ on \ interval:$  In this case, $e_{v}$ must be 1, and $A\cup B$ 
is disconnected in the 
sense of the previous section.  The
pushforward of $\alpha$ along $\St{A}{B}$ vanishes by 
\refP{disconnected}.
\vskip 5pt
\noindent
$v \ free, \ e_{v}=1:$  This means that there are two more edges 
connecting $v$ to vertices in $A\cup B$.  But this is precisely 
the case of \refL{case1hidden} so the pushforward is zero.
\vskip 5pt
\noindent
$v \ free, \ e_{v}=2:$  Now there is another edge emanating from 
$v$ and ending in a vertex in $A\cup B$.  If it ends in an 
interval vertex, $A\cup B$ is again disconnected.  
If it ends in a free 
vertex $v'$, then there are two more edges emanating from this 
vertex.  The situation when both of those end in vertices in $A\cup B$ is 
precisely the setting of \refL{case2hidden}.  If they both end 
in vertices not in $A\cup B$, then $A\cup B$ is disconnected.  
If only one of them 
ends in a vertex in $A\cup B$, then $v'$ with its edges forms a picture as 
in \refL{case1hidden}.
\vskip 5pt
\noindent
$v \ free, \ e_{v}=3:$  $A\cup B$ is disconnected.
\vskip 5pt
Since there cannot be more than three edges emanating from a vertex, 
this completes the proof. 
\end{proof}

\subsection{Principal faces}\label{S:PrincipalFaces}

In this situation of exactly two points colliding, there are 
various cases to consider, depending on whether
the vertices in the diagram $D$ are on the interval 
or free and on whether 
they are connected by a chord or an edge.

An essential difference from the previous arguments will be that 
not all the pushforwards will vanish individually, but rather we will 
have to consider combinations of integrals for various principal faces.  
The combinations are 
determined by the $STU$ and $IHX$ relations.  We will also now 
need to pay attention to the labelings of the diagrams.

Recall from \refT{Thurston} that, given a weight system $W$, the 
expression claimed to vanish is the derivative
\begin{equation}\label{E:sum}
dT(W)=\frac{1}{(2n)!}\sum_{D\in TD_n}W(D)(dI(D,K)+dM_DI(D_1,K)),
\end{equation}
where $dI(D,K)$ is the sum of the pushforwards of $\alpha$ to $\K$ 
along the faces $\St{a}{b}$ of $F[k,s; \K, \R^3]$.  

Now consider the three summands determined by labeled diagrams differing only as 
in the $STU$ relation of Figure \ref{F:STU} in \S\ref{S:Trivalent}.  For simplicity, 
set the labelings as $i=1$ and $j=2$ for now, and we will 
mention more general labelings later.

The first diagram is associated to $F[k-1,s+1; \K, \R^3]$ and the other two to 
$F[k,s; \K, \R^3]$.  Each of 
those has a principal face where points $p_{1}$ and $p_{2}$, corresponding to 
vertices 1 and 2 in the three diagrams, come together.
Denote these faces by $\Ss$, $\ST$, and $\SU$.  The 
goal is to show that the integrals along the three faces
 have the same value, but with signs as in the 
$STU$ relation.  The sum of the three terms in \eqref{E:sum} will then 
be a multiple of
$$
W(S)-W(T)+W(U),
$$ 
and hence 0 by the $STU$ relation.

Let $\alpha_{S}$, $\alpha_{T}$, and $\alpha_{U}$ be the three forms 
which are integrated.  It is clear that
$$
\int\limits_{\ST}\alpha_{T}=-\int\limits_{\SU}\alpha_{U}
$$
since $\ST$ and $\SU$ are diffeomorphic; all the 
maps $h_{ij}$ for those faces are identical, so that 
$\alpha_{T}=\alpha_{U}$; but the orientations on $\ST$ and 
$\SU$ are different since the two labels, 1 and 2, are switched.  
The two integrals are thus the same, but with opposite signs 
since the two vertices are not connected by a chord (had there been a 
chord between them, the orientation would change, but the 
map $h_{12}$ would also become its negative, cancelling the negative 
coming from the orientaton).

\begin{rem} By writing an integral over a face like $\ST$ and $\SU$ 
we mean that the integral is taken over the fiber of that face over 
$\K$.
\end{rem}
It remains to show that, for example,
$$
\int\limits_{\Ss}\alpha_{S}=-\int\limits_{\ST}\alpha_{T}.
$$
But $\Ss$ is a face of $F[k-1,s+1; \K, \R^3]$, while $\ST$ is a face of
$F[k,s; \K, \R^3]$.  There is also an extra map $h_{12}$ on $\Ss$, coming from the 
edge connecting vertices 1 and 2.  However, the remaining maps 
$h_{ij}$ are the 
same on the two faces.

To see what the extra map $h_{12}$ is, note that a neighborhood of $\Ss$ is parametrized by
\begin{equation}\label{E:StratumParameters-Ss}
(x_{1}, u, r\, ;\, x_{2}, \ldots, x_{k+s-1}),\ \ \ u\in S^2,\ r\geq 
0,
\end{equation}
with other conditions imposed on the $x_{i}$ as in 
\eqref{E:StratumConditions-ab}.
The two configuration points corresponding to diagram vertices 1 and 
2 are
$$p_{1}=K(x_{1}), \ \ \ p_{2}=K(x_{1})+ru,
$$
so that the extra map on $\Ss$ is simply given by
$
h_{12}=u.
$
It follows that $\alpha_{S}$ and $\alpha_{T}$ differ only in that 
$\alpha_{S}$ contains one more factor than $\alpha_{T}$, the $(m-1)$-form 
$(h_{12})^{*}\omega_{12}$.  In short,
$$
\alpha_{S}=(h_{12})^{*}\omega_{12}\cdot \alpha_{T}.
$$
But since $h_{12}$ is simply the identity on $S^2$, 
$(h_{12})^{*}\omega_{12}$ may be identified with $\omega_{12}$.  
More concisely, using 
Fubini's Theorem, the pushforward of $\alpha_{S}$ to $\K$ can be rewritten as
\begin{equation}\label{E:fubini}
\int\limits_{\Ss}\alpha_{S}=
\pm\int\limits_{u\in S^{2}}u^{*}\omega_{12}
\int\limits_{x_{i}}\alpha_{T}
=\pm\int\limits_{S^{2}}\omega_{12}\int\limits_{\ST}\alpha_{T}
=\pm\int\limits_{\ST}\alpha_{T}
\end{equation}
The indeterminacy in sign comes from the possibly different orientations 
of $\Ss$ and $\ST$.  Seeing that these are in fact 
opposite is straightforward and we leave this to the reader.

To complete the argument, the parameters for the remaining points in 
$\Ss$ and $\ST$ should be mentioned, although they are of no 
consequence.  The outward normal vector giving orientation on these 
principal faces should always be added to the basis of the tangent 
space as the last vector.  The orientations on $\Ss$ and $\ST$ will 
still come out to be 
different as in the above.  This also takes care of an arbitrary 
labeling of the two vertices considered:  Switching any two labels 
affects the orientations on $\Ss$ and $\ST$ by switching two tangent vectors
or a tangent vector and the outward normal, but the same permutation 
occurs for both spaces.  The orientations on $\Ss$ and $\ST$ will thus 
always be different.

The situation is even simpler for the three diagrams in the $IHX$ 
relation since all three corresponding spaces have the same number of 
points on and off the knot.  It is clear that the integrals of 
$\alpha$ over $\SI$, $\SH$, and $\SX$ will have the same value, and 
the labeling of the vertices will ensure that the signs 
come out as required in the relation.

\vskip 5pt
\noindent
This leaves the 
principal faces determined by the pictures in the figure below.

\begin{figure}[h]
\begin{center}
\input{RemainingPrincipal2.pstex_t}
\caption{Remaining principal faces}
\end{center}
\end{figure}
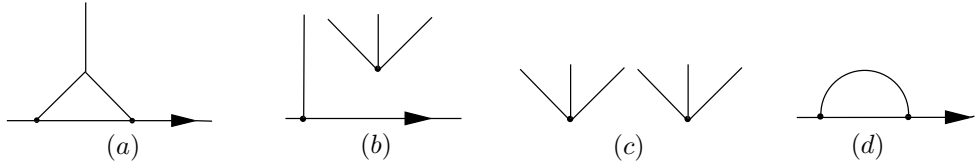
When the two vertices on the interval in figure $(a)$ come together, the 
two maps corresponding to their two edges become the same.  Thus 
$h_{D}$ on this stratum factors through $(S^2)^{|e|-1}$ and the 
pullback of $\omega$ is thus necessarily 0.
If the maps $h_{ij}$ are determined by edges as in figures $(b)$ and 
$(c)$, the integral vanishes 
by \refP{disconnected}.
Finally remember that we are only considering primitive 
weight systems, namely those that vanish on diagrams that can be 
obtained as connected sums of two other diagrams.  Figure (d) 
represents such a summand, so that we may disregard this principal 
face.

\subsection{Faces at infinity}\label{S:InfinityFaces}

Here we 
are saved by the requirement that all our trivalent diagrams are 
connected, so that the set of trivalent vertices is ``anchored'' to 
the interval by at least two edges (it is an easy combinatorial fact 
that connected trivalent diagrams with free vertices 
cannot have just one edge connecting the interval vertices to the free 
ones).  

Denote by $\mathcal{S}_{\infty}$ the stratum with one or more points on the long
 knot or off it going to infinity. 
\begin{prop}\label{P:InfinityFaces} The pushforward of $\alpha$ to 
$\K$ along $\mathcal{S}_{\infty}$ vanishes.
\end{prop}
\begin{proof}
First recall from the end of \S\ref{S:Coordinates&Forms} that the 
parametrization of faces at infinity is much 
like the one of colliding points, except that $r$ approaches $\infty$ 
rather than 0.  Next observe that \refL{case1hidden} still holds for 
$\mathcal{S}_{\infty}$ since the argument used in the proof did not 
depend on $r$ tending to 0.

Given $D$, let $\mathcal{S}_{\infty}$ 
determined by
some subset $B$ of the free vertices of $D$ escaping to infinity.  
These are connected 
to the rest of $D$ by 
some number of external edges.  As in the proof of \refP{hidden}, let 
$v$ be any 
vertex in $B$ connected to one or more vertices in the complement of 
$B$.
Also let $e_{v}$ 
be the number of the external edges of $v$.
Then we have: 
\vskip 5pt
\noindent
$e_{v}=1:$  There are two more edges emanating from 
$v$ which are connected to vertices in $B$.  The picture 
now
is the one of \refL{case1hidden} and so the pushforward of $\alpha$ 
must be 0.
\vskip 5pt
\noindent
$e_{v}=2:$  We could invoke \refL{case2hidden} to dispense with this 
case, but a more direct argument is that there will now be two maps 
$h_{ij}$ which are identical on $\mathcal{S}_{\infty}$.  This is 
immediate from recalling how $h_{ij}$ extends to strata at infinity, 
in particular from \eqref{E:InfinityDirection2}.  The pullback 
$\alpha=(h_{D})^{*}\omega$ to  $\mathcal{S}_{\infty}$ thus factors
through a space of lower dimension than that of $\omega$, so that 
$\alpha$ must be identically 0.
\vskip 5pt
\noindent
$e_{v}=3:$  Now there are three maps $h_{ij}$ which are identical, and 
the comments of the previous case apply again.
\vskip 5pt
\noindent
This exhausts all the cases since a free vertex has exactly three 
edges emanating from it.  
\end{proof}

\subsection{Anomalous Faces}\label{S:AnomalousFaces}

The only case of a pushforward along a codimension one stratum of $F[k,s; \K, \R^3]$ left to 
consider is that of the anomalous face, with all $k+s$ points coming together
 at the 
same time.  
This is the only instance where the situation differs for classical knots and knots in higher codimension.  We summarize some of what is 
known about the anomalous faces and the correction factor $M_DI(D_1,K)$ for classical knots
in \refP{anomalous}.  For $\K_m$, $m>3$, see \refP{AnomalousVanishes}. The observations and 
constructions presented here 
are
direct extensions of those in \cite{BT}.

Recall that the pushforward along the anomalous face  
vanishes if $D$ has at least one chord (\refC{anomalous-cd}).  However, 
this may not be the case for other trivalent
diagrams, so assume that $D$ has $k$ interval vertices, 
$s$ free vertices, and no chords.
A neighborhood of the anomalous face $\mathcal{S}_{D}$ in $F[k,s; \K, \R^3]$ is 
parametrized by
\begin{equation}\label{E:AnomalousParameters}
(x, u_{1}, \ldots, u_{k}, \ldots, u_{k+s}, r),
\end{equation}
where $u_{1}, \ldots, u_{k}\in I$ and $u_{k+1}, \ldots, 
u_{k+s}\in\R^3$ satisfy \eqref{E:StratumConditions-ab}.  

Since $D$ contains no chords, and all 
points are colliding, $h_{D}$ can only be a product of two types of 
$h_{ij}$ upon restriction to $\mathcal{S}_{D}$, namely
$$
h_{ij}=\frac{u_{j}-K(x)}{|u_{j}-K(x)|} \ \ \ \text{and}\ \ \ 
h_{ij}=\frac{u_{j}-u_{i}}{|u_{j}-u_{i}|}.
$$
There is only one configuration point 
left after collision, so $\mathcal{S}_{D}$ maps to $F[1,0; \K, \R^3]$ by projection on 
$x$.  Notice now that 
the pushforward of $\alpha$ may be thought of as
$$\int\limits_{(x,u_{i})}\alpha=\int\limits_{x}\bigg( \ \ 
\int\limits_{\stackrel{u_{i}}{x \text{ fixed}}}\alpha\bigg).$$
In other words, if we let $p$ be the projection as above, and 
$\widehat{\pi}$ the bundle map from $F[1,0; \K, \R^3]$ to $\K$, then the 
diagram
\begin{equation}\label{E:AnomalousTriangle}
\xymatrix{
\mathcal{S}_{D} \ar[dr]^{p_{*}} \ar[dd]_{(\partial\pi)_{*}} &    \\
&   \G{1}{0}  \ar[dl]^{\widehat{\pi}_{*}}                  \\
\K  &  \\
}
\end{equation} 
commutes.  

Now let $N_{D}$ be a subspace of $I^{k}\times\R^{3s}\times 
S^{2}$ consisting of points
$(w_{1}, \ldots, w_{k}, w_{k+1}, \ldots, w_{k+s}, v)$
satisfying
\begin{flalign}
 & w_{1}, \ldots, w_{k}\in I; & \label{E:MDCondition1}\\
 & w_{k+1}, \ldots, w_{k+s}\in\R^3; & \\
 & w_{1}<\cdots <w_{k}; & \\
  & \text{if $i,j>k$ are vertices in $D$ connected by an edge, then 
   $w_{i}\neq w_{j}$}; & \\
 & \text{if $i\leq k$, $j>k$ are vertices in $D$ connected by an 
edge, then
   $w_{i}v\neq w_{j}$};  & \\
 & \sum_{i=1}^{k+s}|w_{i}|^{2}=1; & \\
 & \sum_{i=1}^{k}w_{i}+\sum_{i=k+1}^{k+s}\langle v,w_{i} \rangle 
 =0; &  \label{E:MDCondition7}
\end{flalign}
We now get a commutative diagram
\begin{equation}\label{E:MDDiagram}
\xymatrix{
\mathcal{S}_{D} \ar[r]^{f} \ar[d]^{p} &  N_{D}  \ar[d]^{g}  \\
F[1,0; \K, \R^3] \ar[r]^{\tau} & S^{2} 
}
\end{equation}
where $\tau$ is the usual tangential map
$x\mapsto K'(x)/|K'(x)|,$
$g$ is the projection $(w_{i},v)\mapsto v,$
and $f$ is given by
$$(x, u_{1}, \ldots, u_{k}, u_{k+1}, \ldots u_{k+s}, 0)\longmapsto
  (u_{1}|K'(x)|, \ldots, u_{k}|K'(x)|, u_{k+1}, \ldots, 
  u_{k+s}, \frac{K'(x)}{|K'(x)|}).$$
Because of the usual conditions the parameters 
\eqref{E:AnomalousParameters} have to satisfy, 
$f$ 
preserves \eqref{E:MDCondition1}--\eqref{E:MDCondition7}.

Let $\widehat{h}_{ij}$ denote maps from $N_{D}$ to $S^{2}$ determined by 
the edges of $D$; if $i$ is on the interval and $j$ is free, or if 
both are free, we get respectively 
$$\widehat{h}_{ij}=\frac{w_{j}-w_{i}v}{|w_{j}-w_{i}v|},
\ \ \ 
\widehat{h}_{ij}=\frac{w_{j}-w_{i}}{|w_{j}-w_{i}|}.$$
Then $h_{ij}$ clearly factor through $N_{D}$, so letting 
$\widehat{h}_{D}=\Pi_{edges}\widehat{h}_{ij}$ and $(\widehat{h}_{D})^{*}\omega 
=\widehat{\alpha}$ gives
$$\alpha =h_{D}^{*}\omega =(\widehat{h}_{D}\circ f)^{*}\omega 
=f^{*}((\widehat{h}_{D})^{*}\omega)=f^{*}\widehat{\alpha}.$$
Bott and Taubes call $\widehat{\alpha}$ ``universal'' since it only 
depends on the diagram $D$ and not on the embedding.

Now use \eqref{E:AnomalousTriangle} so that 
$(\partial\pi)_{*}\alpha 
=\widehat{p}_{*}(p_{*}(f^{*}\widehat{\alpha})).$
Further, note that $N_{D}$ was defined precisely so that the diagram in 
\eqref{E:MDDiagram} is 
in fact a pullback square, so we can actually 
obtain $(\partial\pi)_{*}\alpha$ by integrating $\widehat{\alpha}$ 
to $S^2$, pulling the result back via $\tau$ and then integrating 
to $\K$.  In short, the complete diagram we are using is
\begin{equation}
\xymatrix{
      &   S_{D} \ar[rr]^{f} \ar[dd]^{p} \ar[dl]_{\partial\pi} &  & 
N_{D} 
      \ar[r]^{\widehat{h}_{D}} \ar[dd]^{g} &  (S^2)^{|e|} \\
  \K &  & &  &  & \\
      &  F[1,0; \K, \R^3] \ar[ul]^{\widehat{p}} \ar[rr]^{\tau} &   & S^2  &
}
\end{equation}
with
$$(\partial\pi)_{*}\alpha 
=\widehat{p}_{*}(\tau^{*}(g_{*}\widehat{\alpha}), \ \ \ 
\widehat{\alpha}=(\widehat{h}_{D})^{*}\omega.$$

Since $N_{D}$ has dimension $k+3s$, the dimension of its fiber over 
$S^2$ is $k+3s-2$.  On the other hand, $\widehat{\alpha}$ is a 
$(k+3s)$-form as $\omega$ is.  It follows that the integration 
$g_{*}\widehat{\alpha}$ produces a $2$-form on $S^{2}$ 
(remember that $2n=k+s$ is the total number of vertices in $D$).  But this form 
 must be rotationally invariant
 since $\omega$ is.  The only such 
2-forms on $S^{2}$ are constant multiples of the standard unit volume form, 
which we denote by $\omega_{12}$.  

We can now summarize the 
situation with the anomalous faces in the following
\begin{prop}\label{P:anomalous}  For a diagram $D$ with chords, 
the pushforward of $\alpha$ to $\K$ along the anomalous face is zero (\refC{anomalous-cd}).  For $D$ with no chords, the pushforward is
$$\mu_D\int\limits_{F[1,0; \K, \R^3]}\tau^{*}\omega_{12}$$
where 
$\omega_{12}$ is the unit volume form on $S^{2}$, $\tau$ is the 
tangential map from $F[1,0; \K, \R^3]$ to $S^{2}$, and $\mu_{D}$ is a real 
number which depends on $D$.
\end{prop}

Note that this finally proves \refT{Thurston}.

\vskip 5pt
\noindent
To define the correction term $M_DI(D_1,K)$, we now look for a space with 
boundary $F[1,0; \K, \R^3]$ and a map to $S^{2}$ which restricts to the 
tangential map $\tau$ on that boundary.  The answer can be traced back 
to \S\ref{S:LinkingNumber} and is clearly the space $F[2,0; \K, \R^3]$ with 
the normalized difference of the two points, $p_{1}$ and $p_{2}$ giving the map to the 
sphere.  Note that the boundary of $F[2,0; \K, \R^3]$ has two diffeomorphic components, depending on 
which order $p_{1}$ and $p_{2}$ appear on the interval.

\begin{definition}\label{D:Correction} Letting  $M_D=\mu_D/2$ and defining $I(D_1,K)$ as before by
$$
I(D_1,K)=\int\limits_{F[2,0; \K, \R^3]}\left( \frac{p_{2}-p
_{1}}{|p_{2}-p_{1}|} \right)^{*}\omega_{12}
=\int\limits_{F[2,0; \K, \R^3]}h_{12}^{*}\omega_{12}
$$
gives the correction term from \refT{Thurston}.
\end{definition}

This ``anomalous term," which has been conjectured to be zero, has been shown to vanish in even degrees \cite{Alt} and in degrees 3 and 5 \cite{Th, Poir}.  This conjecture has recently been reformulated by S.-W. Yang and C.-H. Yu \cite{YY} in terms of the computation of the homology of the trivalent graph complex ($\Td$ can be turned into a complex via contraction of edges as the boundary operator).
As mentioned in the introduction, the vanishing of the anomaly is intimately related to the question of equivalence of Bott-Taubes integrals and the Kontsevich Integral \cite{Kont}.  Namely, Poirier was able to show that the two are indeed the same up to vanishing of the anomalous term \cite{Poir}.  Using Poirier's work, Lescop \cite{L} additionally proved that the computation of the anomaly only has to be carried out on a certain smaller subclass of diagrams.

\section{Universal Finite Type Invariant}\label{S:Universal}

One of the most striking features of Chern-Simons perturbation theory 
is its relation to finite type knot invariants.  The fact that 
configuration space integrals can be used to construct the universal 
finite type knot invariant has been known for some time \cite{Alt, 
Th}, and we now provide the details of the proof of this fact.  This can be in some sense viewed as completion of 
Thurston's work in \cite{Th}.  More details about finite type 
invariants can be found in \cite{BN, BN2}.

\vskip 5pt
\noindent  
 A \emph{singular knot} 
is a knot as before except for a finite number 
of double points.  The tangent vectors at the double points are required 
to be independent.  A knot with $n$ such
self-intersections is called \emph{$n$-singular}.

Any knot invariant $V$ can be extended to singular knots via a 
repeated use of the \emph{Vassiliev skein relation} pictured in 
Figure \ref{F:VassilievSkein}.

\begin{figure}[h]
\begin{center}
\input{skeinrelation.pstex_t}
\caption{Vassiliev skein relation}\label{F:VassilievSkein}
\end{center}
\end{figure}
The drawings of the knot projections are meant to indicate that the 
three knots only differ locally in one crossing.  The two knots on the 
right side of the equality are called the \emph{resolutions} of a 
singularity.  A singular knot with $n$ singularities thus produces 
$2^{n}$ resolutions, and the order in which the sigularities are 
resolved does not matter due to the sign conventions.

\begin{definition}
\label{D:FiniteTypeInvariant}
$V$ is a \emph{(finite, or Vassiliev) type $n$ invariant} if it 
vanishes identically on singular knots with $n+1$ self-intersections.
\end{definition}

Let $\V$ be the collection of all finite type invariants and let 
$\V_{n}$ be the type $n$ part of $\V$.  An immediate consequence of 
\refD{FiniteTypeInvariant} is that $\V_{n}$ contains $\V_{n-1}$.
It is also clear from the definition that $\V_{0}$
consists only of constant functions on $\K$.

Suppose now that $K_{1}$ and $K_{2}$ are $n$-singular knots with 
singularities in the same place, by which we mean that the points 
on $I$ which the immersions identify in pairs appear in the same order 
on $I$.  It is 
clear that $K_{2}$ may be obtained from $K_{1}$ by a sequence of 
crossing changes.
But if $V$ is type $n$, it follows that 
$
V(K_{1})=V(K_{2}),
$
since the difference of the value of $V$ on $K_{1}$ and its value on 
the same $n$-singular knot as $K_{1}$ but with a crossing changed is 
precisely the value of $V$ on an $(n+1)$-singular knot according to 
the Vassiliev skein relation.  Since $V$ is type $n$, it must 
vanish on such a knot by definition.  We thus note that
\vskip 4pt
\noindent
\emph {The value of a type $n$ invariant on an $n$-singular knot only 
depends on the placement of its singularities.}
\vskip 4pt
\noindent
With this observation, it is clear that any type 1 invariant also 
must be the constant function on $\K$, with the exception of the 
framing number of one considers framed knots.
It can also be shown that, up to framing,
there is a unique nontrivial type 2 invariant.  In fact, 
one interpretation of \refT{Universality} below is that the values 
of all finite type invariants on all knots can be computed 
inductively.  In practice, however, such computations are quite 
complicated \cite{BN}.

Getting back to our observation, the value of a type $n$ invariant 
$V$ on an $n$-singular knot thus only depends on what can schematically be 
represented as an interval with $2n$ paired-off points on it, i.e. a 
\emph{chord diagram}.  The 
pairs 
 serve as prescriptions for where the singularities on the knot should 
occur.

Let $CD_{n}$ be the set of all chord diagrams with $n$ chords.
We thus have that a type $n$ invariant $V$ determines a function on 
$CD_{n}$.  More precisely, if $D$ is an element of $CD_{n}$, and if 
$K_{D}$ is any $n$-singular knot with singularities as prescribed by 
$D$, we have a map
\begin{equation}\label{E:EasyMap}
\V_{n}\longrightarrow \{f\colon \R[CD_{n}]\to\R\}
\end{equation}
given by
\begin{equation}\label{E:EasyDirection}
f(D)=V(K_{D})
\end{equation}
and extending linearly.
Note that the kernel is by definition $\V_{n-1}$.

However, not every function on chord diagrams arises in this way:  
Suppose $V$ is evaluated on a knot with $n-1$ singularities and 
suppose a strand is moved in a circle around one of them.
Along the 
way, the strand will
introduce a singularity each time it passes through another strand emanating from 
the original singularity.  Thus four new $n$-singular knots are 
created in the process.  Since the strand is back where it started, the sum
 of the values of $V$ on the four $n$-singular knots 
should be zero.  But each of the four knots corresponds to a chord 
diagram.  It follows that if $f$ is in the image of the map
(\ref{E:EasyMap}) then the sum of its values on those four 
diagrams is 0.  This is known as the $4T$ (\emph{four-term}) relation, 
which upon unravelling the four knots into chord diagrams with 
appropriate signs appear in Figure \ref{F:4T}.  The diagrams differ 
only in chords indicated; there may be more 
chords with their endpoints on the dotted segments, but they are the same for all four 
diagrams.

\begin{figure}[h]
\begin{center}
\input{4T.pstex_t}
\caption{$4T$ relation}\label{F:4T}
\end{center}
\end{figure}
\begin{rem}
In most literature on finite type theory, one more 
relation besides the $4T$ is imposed.  This $1T$ (\emph{one-term) 
relation} sets the value of any weight system on a chord diagram with 
an \emph{isolated chord}, i.e. a chord not intersected by any other 
chords, to be zero.  
However, if one works with framed knots, $1T$ relation cannot be imposed.  This is because 
the two resolutions of a singularity coming from an isolated chord 
are \emph{not framed isotopic}.
The consequence of having to consider framed knots is simply that the 
number of finite type invariants is somewhat larger.  One now
 gets a genuine type 1 invariant, 
the framing number, whose $n$th power is the additional type $n$ 
invariant.  The anomalous term from \refD{Correction} also has to be modified slightly \cite{Poir, L2} and varies with the framing.
\end{rem}

Let 
$\cd{n}=\R[CD_{n}]/4T$ and recall that 
$\td{n}$ is the real vector space generated by trivalent diagrams modulo
 the $STU$ relation.
The following important theorem, due to Bar-Natan, gives the 
 connection between configuration space integrals and 
finite type invariants:
\begin{thm}[\cite{BN}, Theorem 6]\label{T:Chord=Trivalent}
$\cd{n}$ and $\td{n}$ are isomorphic for all positive $n$, so that every weight system 
$W\in\W_{n}$ extends uniquely from $\cd{n}$ to $\td{n}$.
\end{thm}
The idea of the proof is to construct the map giving an isomorphism 
inductively, noting that, in the base case of one free and three 
interval vertices, the $4T$ relation is the difference of two $STU$ 
relations.

Now we can think of $\W_{n}$, the space of weight systems of 
degree $n$, as those $f$ in above which vanish on the $4T$ relation.  
It turns out that $\W_{n}$ is all there is to the image of the map in 
\eqref{E:EasyMap}.  Its inverse is given by $T(W)$:
\begin{thm}\label{T:Universality} $T(W)$ is a type $n$ knot 
invariant.  Further, it gives an isomorphism between $\W_{n}$ and 
$\V_{n}/\V_{n-1}$.
\end{thm}
This theorem was first proved by Altschuler and Freidel \cite{Alt} in a somewhat different form.  Before proving it, we need the following

\begin{lemma}\label{L:UniversalityLemma}

Let $D'\in CD_{n}$ be a labeled chord diagram and let $K_{D'}$ be 
any $n$-singular knot, $n>1$, with singularities as prescribed by 
$D'$.  Also 
let $K_{D',S}$ be the resolutions of $K_{D'}$ determined by nonempty 
subsets $S$ of $\{1, \ldots,n\}$ where the $i$th singularity is 
resolved positively (the first resolution in the Vassiliev skein 
relation) if $i\in S$.  Then, given $\delta>0$, there are isotopies 
of
the $K_{D',S}$ to knots $K_{D',S}'$ such that

\begin{align}
&\sum_{K_{D',S}'}M_DI(D_1,K)=0, 
  \ \ \ \ \ \ \ \ \ \ \ \ \ \text{ for all labeled $D\in TD_{n}$,\ $n>2$} 
  \label{E:UniversalityLemma1} \\ 
&\Big|\sum_{K_{D',S}'}M_DI(D_1,K)\Big|<\delta, 
  \ \ \ \ \ \ \ \ \ \ \ \text{ for all labeled $D\in TD_{2}$,} 
  \label{E:UniversalityLemma2} \\ 
&\Big|\sum_{K_{D',S}'}I(D, K_{D',S}')\Big|<\delta, 
  \ \ \ \ \ \ \ \ \ \ \ \ \text{ for all labeled $D\in TD_{n}$,\ $D\neq D'$}. 
  \label{E:UniversalityLemma3}\\ 
&\Big| \sum_{K_{D',S}'}I(D', K_{D',S}')-1\Big|< \delta. 
  \label{E:UniversalityLemma4} 
\end{align}

\end{lemma}

\begin{proof}  The proof is lengthy but not difficult.
We will simply use the fact that we may 
choose the 
resolutions so that they differ only inside of $n$ balls in $\R^{3}$ 
of arbitrarily small radius.

To show \eqref{E:UniversalityLemma1}, 
consider the sum of anomalous integrals $M_DI(D_1,K)$ over all resolutions $K_{D',S}$ 
(with signs 
according to the skein relation).  The resolutions are isotopic to 
knots which are the same outside of some disjoint balls 
$B_{i}$ in $\R^{3}$, $1\leq i\leq n$.  
In other words, if $a_{1}, b_{1}, a_{2}, 
b_{2}, \ldots, a_{n}, b_{n}$ are points on $I$ which make up the $n$
singularities, the isotopy is given be reparametrizing the interval 
so that smaller 
neighborhoods $N_{i}$ of each $b_{i}$ are 
embedded as either the overstrand or the understrand when making the 
resolutions.  All resulting resolutions $K_{D',S}'$ can then be 
paired off, with 
opposite signs, into knots differing only in $K_{D}'(N_{i})$.  We 
denote the embedding of 
$N_{i}$ as the overstrand by 
$K_{D}'(N_{i})^{+}$ and by $K_{D}'(N_{i})^{-}$ otherwise.
Each $B_{i}$ can be chosen so that it contains both 
$K_{D}'(N_{i})^{+}$ and 
$K_{D}'(N_{i})^{-}$ and so that it is disjoint from all other 
$B_{j}$.  Notice that by arranging a suitable isotopy,  
$B_i$ can be made to have arbitrarily small radius and contain 
no other parts of the knot besides the arcs $K_{D}'(N_{i})^{+}$ and 
$K_{D}'(N_{i})^{-}$.

We may now break up $M_DI(D_1,K)$ for each resolution
into a sum of 
integrals over various neighborhoods of the configuration space 
$\FM{2}{K_{D',S}'}$.  These are determined by whether $p_{1}$ and $p_{2}$ 
are in some $B_{i}$ and $B_{j}$.  If they are 
both 
outside of $B_{i}$ for some $i$, then the resolutions can be 
paired off so that 
each pair differs only inside of $B_i$.  The two integrals in each 
pair, taken over such a 
neighborhood, are identical and appear with opposite signs due to 
the skein relation, so they cancel. 

If one or both of $p_{1}$ and $p_{2}$ are inside some $B_i$ and 
$B_j$, a similar 
situation occurs.  Since $n>2$, we can always pair off the 
resolutions so that each pair differs only inside of some third ball 
$B_k$ which contains neither
$p_{1}$ nor $p_{2}$.    The integrals again cancel after they are 
paired according to how they differ inside $B_{k}$.
%



The argument for \eqref{E:UniversalityLemma2} is more complicated  
since there are now exactly two balls $B_{1}$ and 
$B_{2}$.  As $p_{1}$ and $p_{2}$ may fall into both of them, 
there are various neighborhoods (depending 
on whether the points are on the overstrands 
or the understrands)
of $\FM{2}{K_{D',S}'}$ for which the integrals do not cancel.
However, we will show that they can be paired so that their difference is 
arbitrarily small.

For example, pick the resolution $K_{D',\{1,2\}}'$ and a neighborhood of 
$\FM{2}{K_{D',\{1,2\}}'}$ 
where $p_{1}$ is on the arc $K_{D}'(N_{1})^{+}$ and 
the other point, which we call $p_{2}^{+}$, is on 
$K_{D}'(N_{2})^{+}$.  
Then there is another resolution 
$K_{D',\{1\}}'$ and a neighboorhood of $\FM{2}{K_{D',\{1\}}'}$ such that 
$p_{1}$ is again on the arc $K_{D}'(N_{1})^{+}$ while the other 
point, which we now call
$p_{2}^{-}$, is on $K_{D}'(N_{2})^{-}$.

The domains of the integrals $M_DI(D_1,K)$ over this neighborhood are 
$I\times I$ for every resolution since each arc containing a configuration 
point is 
diffeomorphic to $I$.  Further, the integrals appear with opposite 
signs again due to the skein relation, so that we are now studying the difference
\begin{equation}\label{E:IntegralDifference}
\int\limits_{I\times 
I}\left(\Big(\frac{p_{2}^{+}-p_{1}}{|p_{2}^{+}-p_{1}|}\Big)^{*}\omega_{12}-
\Big(\frac{p_{2}^{-}-p_{1}}{|p_{2}^{-}-p_{1}|}\Big)^{*}\omega_{12}\right).
\end{equation}



The difference of the 2-forms in the integral can be written as
$$
\frac{1}{|p_{2}^{+}-p_{1}|^{3}}\, det(\dot{p_{1}}, \dot{p_{2}}^{+}, 
p_{2}^{+}-p_{1})-
\frac{1}{|p_{2}^{-}-p_{1}|^{3}}\, det(\dot{p_{1}}, \dot{p_{2}}^{-}, 
p_{2}^{-}-p_{1})
$$
where the derivatives are taken with respect to the knot parameter 
$t$.  
Now rewrite the above as
\begin{gather}
\frac{1}{|p_{2}^{+}-p_{1}|^{3}}  
\Big(det(\dot{p_{1}}, \dot{p_{2}}^{+}, 
p_{2}^{+}-p_{1})-det(\dot{p_{1}}, \dot{p_{2}}^{-}, 
p_{2}^{-}-p_{1})\Big) \notag \\
+ 
\left(\frac{1}{|p_{2}^{+}-p_{1}|^{3}}-\frac{1}{|p_{2}^{-}-p_{1}|^{3}}   
\right)
det(\dot{p_{1}}, \dot{p_{2}}^{-}, 
p_{2}^{-}-p_{1}). \label{E:LinearDeterminant}
\end{gather}
To show that the second term can be made small by isotoping the 
resolutions, we first need to 
bound the derivatives while the isotopies are performed.  To do this, 
choose smooth bump functions $f(t)$ and $g(t)$ and scale them by
$$
f_{\rho}(t)=\rho^{2}f\left(\frac{t}{\rho}\right),\ \ \ \ 
g_{\rho}(t)=\rho^{2}g\left(\frac{t}{\rho}\right).
$$
The resolutions can be chosen so that the parametrizations for 
the 
points are
\begin{equation*}
p_{2}^{+}=(t,0,f_{\rho}(t)),  \ \ \  
p_{2}^{-}=(t,0,-f_{\rho}(t)), \ \ \ 
p_{1}=(t,0,g_{\rho}(t))
\end{equation*} 
The derivatives $\dot{p_{2}}^{+}, \dot{p_{2}}^{-}$, and 
$\dot{p_{1}}$ are now all bounded so that in particular there is a 
bound $M$ on the determinant in the second term of 
\eqref{E:LinearDeterminant}:
$$
det(\dot{p_{1}}, \dot{p_{2}}^{-}, p_{2}^{-}-p_{1})<M.
$$
But  
the resolutions can also be changed by isotopies so that the distance from $p_{2}^{+}$ to 
$p_{2}^{-}$ is small compared to the distance between $p_{1}$ and 
either of those points.  Namely, we can arrange  for 
$$
\frac{1}{|p_{2}^{+}-p_{1}|^{3}}-\frac{1}{|p_{2}^{-}-p_{1}|^{3}}<
\frac{\epsilon}{M}.
$$
This isotopy shrinks the balls $B_{1}$ and $B_{2}$.

For the first term in \eqref{E:LinearDeterminant}, we can use the 
linearity of the determinant to rewrite it as
$$
det(\dot{p_{1}}, \dot{p_{2}}^{+}-\dot{p_{2}}^{-}, p_{2}^{+}-p_{1})+
det(\dot{p_{1}}, \dot{p_{2}}^{-}, p_{2}^{+}-p_{2}^{-}).
$$
Here we have that
\begin{gather*}
\text{$\dot{p_{1}},
\dot{p_{2}}^{+},$ and  
$\dot{p_{2}}^{-}$ are bounded because of 
$f_{\delta}(t)$ and $g_{\delta}(t)$,} \\
\text{$\frac{1}{|p_{2}^{+}-p_{1}|^{3}}$ and $p_{2}^{+}-p_{1}$ are 
bounded
 because of the isotopy, and} \\
 \text{$p_{2}^{+}-p_{2}^{-}$ is small because of the isotopy.}
\end{gather*}
That leaves the difference 
$$
\dot{p_{2}}^{+}-\dot{p_{2}}^{-}
=\left(1,0,\rho f'\left(\frac{t}{\rho}\right)\right)-
\left(1,0,-\rho f'\left(\frac{t}{\rho}\right)\right).
$$
But this difference can be made small through a choice of $\rho$.  
Hence the whole first term in \eqref{E:LinearDeterminant} can be made 
smaller than $\epsilon$.  Statement \eqref{E:UniversalityLemma2} 
follows by choosing $\delta=2\epsilon$.

As for \eqref{E:UniversalityLemma3}, assume first $D$ is a chord 
diagram different from $D'$.  Again the resolutions $K_{D',S}$ of an 
$n$-singular knot are isotopic to some
 $K_{D',S}'$ which differ only inside of disjoint balls $B_{i}$, 
 $1\leq i\leq n$, of small radii.

We can break up each configuration space into neighborhoods as 
before, 
so that if all the points $p_{k}$ are outside of some $B_{i}$, the 
integrals cancel in pairs.  Otherwise, let $U$ be a neighborhood where each 
$B_{i}$ contains at least one point.  We then have the following 
cases:

\begin{enumerate}
\item  If not all the points are in 
the balls, then there exists a $p_{k}$ in some $B_{i}$ which is connected 
to a point $p_{k}'$ outside of $B_{i}$.  Assume
$p_{k'}$ is either inside another ball 
$B_{j}$ or ``bounded away'' from $B_{i}$ by another ball 
as in Figure \ref{Fi:PossibleCases}. 

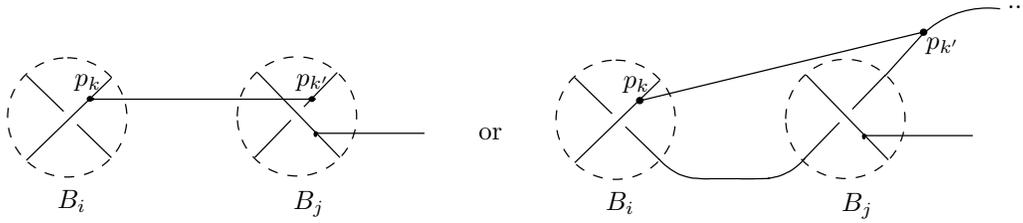
\begin{figure}[h]
\begin{center}
\input{Universality3.pstex_t}
\caption{Case (1) for $D$ a chord diagram different from $D'$.  A line segment between two points indicates that they are 
related by a map $h_{p_kp_{k'}}$.  We will say the two points are ``connected''.
}
\label{Fi:PossibleCases}
\end{center}
\end{figure}

\noindent
For each of such neighborhoods, 
 we can pair off the integrals which differ only inside of 
$B_{i}$.  Using the same arguments as in the proof of 
\eqref{E:UniversalityLemma2}, the differences 
\begin{equation}\label{E:Differences}
\left(\frac{p_{k}^{+}-p_{k'}}{|p_{k}^{+}-p_{k'}|}\right)^{*}\omega_{k'k^{+}}-
\left(\frac{p_{k}^{-}-p_{k'}}{|p_{k}^{-}-p_{k'}|}\right)^{*}\omega_{k'k^{-}}
\end{equation}
can be made small, and we can consequently 
arrange 
\begin{equation}\label{E:LemmaIntegralDifference}
\Big| \int\limits_{U}\prod_{\substack{\text{chords $ij$} \\ ij\neq k'k}}
\left( \frac{p_{j}-p_{i}}{|p_{j}-p_{i}|}  \right)^{*}\omega_{ij}
\int\limits_{I\times I}\left(  
\left(\frac{p_{k}^{+}-p_{k'}}{|p_{k}^{+}-p_{k'}|}\right)^{*}\omega_{k'k^{+}}-
\left(\frac{p_{k}^{-}-p_{k'}}{|p_{k}^{-}-p_{k'}|}\right)^{*}\omega_{k'k^{-}}  
\right)\Big|<\epsilon
\end{equation}
for each pair of integrals and for any $\epsilon>0$.  Since there are 
$2^{n}$ resolutions which have 
been paired off, we have
$$
\Big|\sum_{K_{D',S}'}I(D, K_{D',S}')\Big|<2^{n-1}\epsilon.
$$

\item Assume again $p_{k}$ is the only point in $B_{i}$, but it is 
connected to a point $p_{k}'$ near or inside $B_{i}$, and the two points are on the same 
strand as indicated in the left picture of Figure \ref{Fi:ImpossibleCases}.  Assume also that either the same situation occurs in every other ball or the ball
 contains two connected points lying on different strands.  
 Otherwise, we could refer to the situations which have been taken care of 
 already.
%


For this, isotope 
the resolutions so that the difference between the overcrossing and 
the undercrossing in $B_{i}$ is contained in a ball of smaller radius 
than that of $B_{i}$.  The differences \eqref{E:Differences} can 
again be made smaller than any positive number by choosing an 
appropriate isotopy, and the difference in integrals differing in 
$B_{i}$ can thus be made to satisfy \eqref{E:LemmaIntegralDifference}.

\item  Suppose each $B_{i}$ contains a 
point $p_{k}$ which is connected to a point $p_{k}'$ near or inside $B_{i}$. 
The two points are on different strands, as indicated in the right picture of Figure \ref{Fi:ImpossibleCases}.
However, this situation occurs only when $D=D'$.
\vskip 11pt

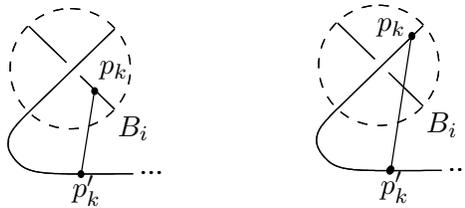
\begin{figure}[h]
\begin{center}
\input{Universality4.pstex_t}
\caption{Cases (2) and (3)}{\label{Fi:ImpossibleCases}}
\end{center}
\end{figure}

\item The the last remaining case is that of all the points in all the balls. Again, because $D$ is different from $D'$, there must be at least two points in different balls which are connected, and we may refer to the first picture in Figure \ref{Fi:PossibleCases}.  However, before being able to argue as in that case, it has to be noted that there may be more than one difference \eqref{E:Differences} 
for each 
pair of integrals since there may be more than one point $p_{k}$ in 
$B_{i}$.  Consequently, proceeding as we did in that situation will only work if $B_i$ does not contain another pair of connected points on different strands.
If this were the case, there would be no way to 
make the 
difference of pullbacks small for this map.  But we may assume this is not the case since the same cannot happen in $B_{j}$ (if 
it did, 
then there would have to be at least three points in each of $B_{i}$ 
and $B_{j}$ which necessarily means that there is a third ball $B_k$ 
with no 
points at all; integrals could be paired according to how they differ 
in $B_k$ and canceled).  We could 
then pair the integrals which differ only inside of $B_{j}$ and 
proceed as before.

\end{enumerate}

\begin{rem}  For most neighborhoods of $\FM{2n}{K_{D',S}'}$ with 
points in 
the $B_{i}$, the integrals in fact cancel.
The only cases when this cannot be done is when $n$ of the 
configuration points are on the arcs $K_{D',S}'(N_{i})^{+}$ for each 
$i$, since these are the only arcs
in which the resolutions differ .  If any one of these arcs is 
free of configuration points, we can 
couple the resolutions which differ only in the ball containing 
that arc and cancel the integrals.  Otherwise, 
\eqref{E:LemmaIntegralDifference} is the best that can be done.

\end{rem}

Before showing how the sum of the integrals can be made small for 
trivalent diagrams that are not chord diagrams, we prove
the last statement of the lemma.  The difference between what was 
argued so far and this case is that $\FM{2n}{K_{D',S}'}$ now has neighborhoods 
as in Case (3) above.  But no amount of ``shrinking'' 
of 
the $B_{i}$ will then make the differences \eqref{E:Differences} 
small.
Changing any of the overcrossings to an undercrossing can be thought of 
as passing the knot through itself, but as the knot does so, the 
map changes significantly, regardless of how small the ball 
$B_{i}$ is.
However, each integral over such a neighborhood 
is the product
$$
\prod_{\text{chords 
$ij$}}\ \int\limits_{S^{2}_{+}\text{ or }S^{2}_{-}}
\left( \frac{p_{j}-p_{i}}{|p_{j}-p_{i}|}  \right)^{*}\omega_{ij}\, 
=\, 
\pm\frac{1}{2^{n}},
$$
where $S^{2}_{+}$ and $S^{2}_{-}$ denote the two hemispheres of 
$S^{2}$.  The sign depends on how many singularities are resolved 
negatively (so that there will be as many maps to the lower 
hemisphere) as well as the labeling of $D'$, which determines the 
orientation of each $\FM{2n}{K_{D',S}'}$.  For example, these always 
combine to give a positive sign if all 
singularities are resolved positively.

Now we add the integrals $I(D', K_{D',S}')$ for all the
resolutions, but also with negative signs if an odd number of 
singularities is 
resolved negatively (from the skein relation).  The result is
 a sum of $2^{n}$ terms, each with the value $1/2^{n}$.  Since all other 
 neighborhoods contribute terms that can be made arbitrarily small, we 
obtain
 $$
 \Big|\sum_{K_{D',S}'}I(D', K_{D',S}')-1\Big|< \delta
 $$
 as desired.

\noindent
To prove the rest of \eqref{E:UniversalityLemma3}, let
 $D$ be a trivalent diagram with $k$ 
interval vertices and $s=2n-k>0$ free ones.  The integrals are now 
taken over the 
fibers of $\G{k}{s}=F[k,s\, ;\, \K,\R^{3}]$ over $\K$.  The fibers 
are different for each resolution $K_{D',S}'$.  However, if we consider the neighborhoods of the fibers with all configuration points on the knot outside 
of at least one of the balls $B_{i}$, the integrals can 
be paired off according to how they differ in $B_{i}$.  They will 
cancel, since the neighborhood chosen is the same for each pair of integrals.

Consider then the neighborhoods in $F[k,s\, ;\, \K_{D',S}',\R^{3}]$ for each $S$ such that some or all of the points on the 
knot are in all of the $B_{i}$.
Since $k<2n$, there has to 
be a $B_{i}$ 
with exactly one point $p_{j}$ on the knot.  If there was no such $B_{i}$,
$D$ would be a chord diagram with $k=2n$.  So $p_{j}$ is 
connected to another point $p_{j'}$ outside of $B_{i}$.
If $p_{j'}$ is on the knot, we are done by the cases considered when 
$D$ was a chord diagram, and this leaves the case when $p_{j'}$ is off the knot.

Pairing off the 
resolutions according to how they differ in $B_{i}$ and then trying to make 
the differences of integrals $I(D, K_{D',S}')$ small will not work now, since 
$p_{j'}$ 
can be anywhere in $\R^{3}$.  But each $B_{i}$ contains at least one point on the knot.  Suppose $a$ of 
those are associated to at least two.
The remaining $n-a$ balls each contain exactly one point on the knot.  
We have thus far argued that unless all of them are connected to 
points off the knot, the integrals can be compared in pairs and we are 
finished.  The 
total number of points accounted for so far is then at least
$2a+(n-a)=n+a.$  This leaves at most $n-a$ points in $\R^{3}$ (the 
total has to be $2n$), and 
since at least $n+a$ knot points are connected to them, it follows 
that there has to be at least one point in $\R^{3}$ connected to 
points in two of the $n-a$
balls as in Figure \ref{Fi:Universality8Figure}.

\begin{figure}[h]
\begin{center}
\input{Universality8.pstex_t}
\caption{}
\label{Fi:Universality8Figure}
\end{center}
\end{figure}

\noindent
For the neighborhood where $p_{i}$ enters $B_{k}$, we can use $B_{k'}$ 
to pair the resolutions, and vice versa.  The neighborhoods are 
diffeomorphic for all 
resolutions, so all that remains is to show that the difference of 
forms which are integrated can be made small.  This can be argued 
exactly the same way as before.  The only slight difference is that a 
point on the knot is replaced by a point in $\R^{3}$.

If $p_{i}$ is neither in $B_{k}$ nor $B_{k'}$, we can 
enlarge the balls as in Figure \ref{Fi:EnlargingBalls}.
Three cases can now be considered:
If $p_{i}$ is in $\widetilde{B}_{k}\!\setminus\! B_{k}$, $B_{k'}$ can be 
used for 
comparing the resolutions.  Similarly for $p_{i}$ in $\widetilde{B}_{k'}
\!\setminus \! B_{k'}$.  If $p_{i}$ is in the complement of 
$\widetilde{B}_{k}\cup \widetilde{B}_{k'}$, either $B_{k}$ or 
$B_{k'}$ can be used.

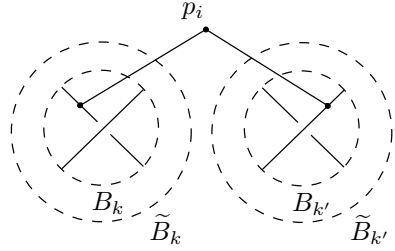
\begin{figure}[h]
\begin{center}
\input{Universality9.pstex_t}
\caption{The remaining neighborhoods for $D$ a trivalent diagram}
\label{Fi:EnlargingBalls}
\end{center}
\end{figure}

\end{proof}

\begin{proof}[Proof of \refT{Universality}]  With the above lemma 
in hand, it is easy to prove that the compositions give the 
identity.  We consider one of them: 
$$
\xymatrix{
\W_{n}\ar[r]^(0.4){T(W)} & \V_{n}/\V_{n-1} \ar@{^{(}->}[r] & \W_{n},
}
$$
where the second map is given by \eqref{E:EasyDirection}.
It suffices to show:  
Given $W\in\W_{n}$ and a chord diagram $D'\in CD_{n}$,
$W(D')=T(W,K_{D'}),$ where

\begin{equation}
T(W,K_{D'})=\frac{1}{(2n)!}\sum_{K_{D',S}}
\sum_{\substack{D\in TD_n \\ D \text{ labeled}}}W(D)(I(D, 
K_{D',S})-M_DI(D_1,K)).
\label{E:EnoughForUniversality}
\end{equation}

One way to prove this is to show 
\begin{equation}\label{E:ModifiedEnough}
\frac{1}{(2n)!}\sum_{K_{D',S}}
\sum_{\substack{D\in TD_n \\ D \text{ labeled}}}D\, (I(D, 
K_{D',S})-M_DI(D_1,K))\ =\ \frac{1}{(2n)!}\Big(
\sum_{\substack{ \text{labelings} \\ 
\text{of $D'$}}}(1+\delta)D'+
\delta\sum_{K_{D',S}}
\sum_{\substack{D\neq D' \\ D\in TD_n \\ D \text{ labeled}}} D \Big),
\end{equation}
where $\delta$ can be made arbitrarily small by isotoping the knots 
$K_{D',S}$.
Since $T(W,K_{D'})$ is an invariant, applying $W$ to both sides of 
\eqref{E:ModifiedEnough} will yield that $\delta$ must be 0.

A simplification can be made by choosing a definite labeling for 
each diagram and then showing 
\begin{equation}\label{E:ModifiedEnough2}
\sum_{K_{D',S}}
\sum_{D\in TD_n}D\, (I(D, 
K_{D',S})-M_DI(D_1,K))\ =\ (1+\delta)D'+
\delta\sum_{K_{D',S}}
\sum_{\substack{D\neq D' \\ D\in TD_n }} D.
\end{equation}
Adding over all labelings and dividing 
by $(2n)!$ will give \eqref{E:ModifiedEnough}.  Moreover, 
\eqref{E:ModifiedEnough2} will follow if, given 
$\delta>0$, there are isotopies of the 
resolutions $K_{D',S}$ such that:
\begin{align}
 & \Big|\sum_{K_{D',S}}(I(D', K_{D',S})-M_{D'}I(D_1,K))-1\Big|<\delta  \notag \\
 & \Big|\sum_{K_{D',S}}(I(D, K_{D',S})-M_DI(D_1,K))\Big|<\delta,  \ \ \ \ \ \ \ \ \ \text{ for all 
$D\neq D'$}. \notag
\end{align}
But this is essentially the statement of \refL{UniversalityLemma} for 
$n>1$.  
Noting that 
in the case $n=1$ the desired statement is trivially true since there 
is only 
one chord diagram with 
one chord, this completes the proof.
\end{proof}

\begin{rem}
The original proof of the isomorphism in \refT{Universality} is due 
to Kontsevich \cite{Kont}.  The role of Bott-Taubes integration is 
played by the 
Kontsevich Integral mentioned in the introduction. 
\end{rem}

\section{Knots in $\R^m$, $m>3$}\label{S:HigherCodimension}

Let $\K_m$ denote the space of knots in $\R^m$, $m>3$.  Bott-Taubes integrals can be defined for $\K_m$ in exactly the same way as described here.  This was first done by Cattaneo, Cotta-Ramusino, and Longoni \cite{Catt} who show that, for a nontrivial weight system $W\in\W_{n}$, 
$T(W)$ produces a nontrivial element of $H^{(m-3)n}(\K_{m})$.  They in fact do more and define a chain complex structure on $\Td$ via contraction of edges to prove

\begin{thm}\label{T:CCL}  There is a chain map $\Td\to \Omega^*(\K)$ which induces 
an injection via Bott-Taubes integration.
\end{thm}
By showing that the diagram complex has nontrivial cohomology in arbitrarily large degrees, one then gets
\begin{cor}
For any $k>0$ and $m>3$, there exists a $l>k$ with $H^l(\K_m)\neq 0$.  
\end{cor}
The main ingredient in the proof of \refT{CCL} is that \refT{Thurston} does not depend on the dimension of 
the Euclidean space.  To show that, one proceeds in exactly the same way as we have in \S\ref{S:Vanishing}.  One only has to be careful with the signs since, depending on whether we are working in $\R^{m}$ for $m$ even or odd, integrals may change sign if the 
configuration points are permuted (corresponding to a permutation of 
the vertices of a diagram) or if the order of those maps is permuted (the 
chords and edges are permuted).  These differences in 
signs then force a definition of two classes of diagrams \cite{Catt, Long2, Long}.  

The odd class (the one defined here) differs from the even class in that one also has to label the edges in the even case, and an appropriate sign has to be introduced in the $STU$ relation.  The $IHX$ relation of \refP{IHX} stays the same (with compatibly labeled edges for the three diagrams), but it is not known if the closure relation now holds.  This means that it is not clear whether the Hopf algebra $\Td$ is commutative in the even case.

With these sign conventions, it is clear that the vanishing and cancelation arguments from \S\ref{S:Vanishing} go through the same way.  One useful observation, however, is that the manifold $N_D$ from \S\ref{S:AnomalousFaces} has dimension $k+ms+m-3$, so that its fiber dimension over $S^
{m-1}$ is $k+ms-2$.  However, $\widehat{\alpha}$ is now a 
$((k+3s)(m-1)/2)$-form.  Then $g_{*}\widehat{\alpha}$ gives a $((m-3)n+2)$-form on $S^{m-1}$.  But this form then
must be 0 if $m>3$ for dimensional reasons.  Thus we get
\begin{prop}\label{P:AnomalousVanishes}
The pushforward of $\alpha$ to $\K_m$, $m>3$, along the anomalous face vanishes for any diagram.
\end{prop}

\end{document}

%% file: LabeledSTU3.pstex_t
\begin{picture}(0,0)%
\includegraphics{LabeledSTU3.pstex}%
\end{picture}%
\setlength{\unitlength}{4144sp}%
\begingroup\makeatletter\ifx\SetFigFont\undefined%
\gdef\SetFigFont#1#2#3#4#5{%
  \reset@font\fontsize{#1}{#2pt}%
  \fontfamily{#3}\fontseries{#4}\fontshape{#5}%
  \selectfont}%
\fi\endgroup%
\begin{picture}(5192,1398)(408,-801)
\put(1880, 29){\makebox(0,0)[lb]{\smash{{\SetFigFont{11}{13.2}{\familydefault}{\mddefault}{\updefault}{\color[rgb]{0,0,0}$=$}%
}}}}
\put(3920, 57){\makebox(0,0)[lb]{\smash{{\SetFigFont{11}{13.2}{\familydefault}{\mddefault}{\updefault}{\color[rgb]{0,0,0}$-$}%
}}}}
\put(5075,-538){\makebox(0,0)[lb]{\smash{{\SetFigFont{12}{14.4}{\familydefault}{\mddefault}{\updefault}{\color[rgb]{0,0,0}$j$}%
}}}}
\put(988,152){\makebox(0,0)[lb]{\smash{{\SetFigFont{12}{14.4}{\familydefault}{\mddefault}{\updefault}{\color[rgb]{0,0,0}$j$}%
}}}}
\put(4345,-532){\makebox(0,0)[lb]{\smash{{\SetFigFont{11}{13.2}{\familydefault}{\mddefault}{\updefault}{\color[rgb]{0,0,0}$i$}%
}}}}
\put(3133,-524){\makebox(0,0)[lb]{\smash{{\SetFigFont{12}{14.4}{\familydefault}{\mddefault}{\updefault}{\color[rgb]{0,0,0}$i$}%
}}}}
\put(2534,-518){\makebox(0,0)[lb]{\smash{{\SetFigFont{11}{13.2}{\familydefault}{\mddefault}{\updefault}{\color[rgb]{0,0,0}$j$}%
}}}}
\put(913,-531){\makebox(0,0)[lb]{\smash{{\SetFigFont{11}{13.2}{\familydefault}{\mddefault}{\updefault}{\color[rgb]{0,0,0}$i$}%
}}}}
\put(823,-746){\makebox(0,0)[lb]{\smash{{\SetFigFont{12}{14.4}{\familydefault}{\mddefault}{\updefault}{\color[rgb]{0,0,0}$S$}%
}}}}
\put(4626,-745){\makebox(0,0)[lb]{\smash{{\SetFigFont{12}{14.4}{\familydefault}{\mddefault}{\updefault}{\color[rgb]{0,0,0}$U$}%
}}}}
\put(2755,-752){\makebox(0,0)[lb]{\smash{{\SetFigFont{12}{14.4}{\familydefault}{\mddefault}{\updefault}{\color[rgb]{0,0,0}$T$}%
}}}}
\end{picture}%

%% file: LabeledIHX2.pstex_t
\begin{picture}(0,0)%
\includegraphics{LabeledIHX2.pstex}%
\end{picture}%
\setlength{\unitlength}{4144sp}%
\begingroup\makeatletter\ifx\SetFigFont\undefined%
\gdef\SetFigFont#1#2#3#4#5{%
  \reset@font\fontsize{#1}{#2pt}%
  \fontfamily{#3}\fontseries{#4}\fontshape{#5}%
  \selectfont}%
\fi\endgroup%
\begin{picture}(5078,1709)(257,-1197)
\put(1288,-253){\makebox(0,0)[lb]{\smash{{\SetFigFont{12}{14.4}{\familydefault}{\mddefault}{\updefault}{\color[rgb]{0,0,0}$=$}%
}}}}
\put(3680,-231){\makebox(0,0)[lb]{\smash{{\SetFigFont{12}{14.4}{\familydefault}{\mddefault}{\updefault}{\color[rgb]{0,0,0}$-$}%
}}}}
\put(653,-1141){\makebox(0,0)[lb]{\smash{{\SetFigFont{12}{14.4}{\familydefault}{\mddefault}{\updefault}{\color[rgb]{0,0,0}$I$}%
}}}}
\put(2522,-1148){\makebox(0,0)[lb]{\smash{{\SetFigFont{12}{14.4}{\familydefault}{\mddefault}{\updefault}{\color[rgb]{0,0,0}$H$}%
}}}}
\put(4614,-1133){\makebox(0,0)[lb]{\smash{{\SetFigFont{12}{14.4}{\familydefault}{\mddefault}{\updefault}{\color[rgb]{0,0,0}$X$}%
}}}}
\put(1674,-239){\makebox(0,0)[lb]{\smash{{\SetFigFont{12}{14.4}{\familydefault}{\mddefault}{\updefault}{\color[rgb]{0,0,0}$-$}%
}}}}
\put(5105,-564){\makebox(0,0)[lb]{\smash{{\SetFigFont{12}{14.4}{\familydefault}{\mddefault}{\updefault}{\color[rgb]{0,0,0}$j$}%
}}}}
\put(4181,-575){\makebox(0,0)[lb]{\smash{{\SetFigFont{12}{14.4}{\familydefault}{\mddefault}{\updefault}{\color[rgb]{0,0,0}$i$}%
}}}}
\put(494, 32){\makebox(0,0)[lb]{\smash{{\SetFigFont{12}{14.4}{\familydefault}{\mddefault}{\updefault}{\color[rgb]{0,0,0}$i$}%
}}}}
\put(495,-518){\makebox(0,0)[lb]{\smash{{\SetFigFont{12}{14.4}{\familydefault}{\mddefault}{\updefault}{\color[rgb]{0,0,0}$j$}%
}}}}
\put(2803,-355){\makebox(0,0)[lb]{\smash{{\SetFigFont{12}{14.4}{\familydefault}{\mddefault}{\updefault}{\color[rgb]{0,0,0}$j$}%
}}}}
\put(2276,-366){\makebox(0,0)[lb]{\smash{{\SetFigFont{12}{14.4}{\familydefault}{\mddefault}{\updefault}{\color[rgb]{0,0,0}$i$}%
}}}}
\end{picture}%

%% file: Closure.pstex_t
\begin{picture}(0,0)%
\includegraphics{Closure.pstex}%
\end{picture}%
\setlength{\unitlength}{4144sp}%
\begingroup\makeatletter\ifx\SetFigFont\undefined%
\gdef\SetFigFont#1#2#3#4#5{%
  \reset@font\fontsize{#1}{#2pt}%
  \fontfamily{#3}\fontseries{#4}\fontshape{#5}%
  \selectfont}%
\fi\endgroup%
\begin{picture}(4642,873)(75,-169)
\put(973,-20){\makebox(0,0)[lb]{\smash{\SetFigFont{12}{14.4}{\familydefault}{\mddefault}{\updefault}\special{ps: gsave 0 0 0 setrgbcolor}. . .\special{ps: grestore}}}}
\put(3133,-21){\makebox(0,0)[lb]{\smash{\SetFigFont{12}{14.4}{\familydefault}{\mddefault}{\updefault}\special{ps: gsave 0 0 0 setrgbcolor}. . . \special{ps: grestore}}}}
\put(2256,219){\makebox(0,0)[lb]{\smash{\SetFigFont{12}{14.4}{\familydefault}{\mddefault}{\updefault}\special{ps: gsave 0 0 0 setrgbcolor}$=$\special{ps: grestore}}}}
\end{picture}

%% file: Type2CD.pstex_t
\begin{picture}(0,0)%
\includegraphics{Type2CD.pstex}%
\end{picture}%
\setlength{\unitlength}{4144sp}%
\begingroup\makeatletter\ifx\SetFigFont\undefined%
\gdef\SetFigFont#1#2#3#4#5{%
  \reset@font\fontsize{#1}{#2pt}%
  \fontfamily{#3}\fontseries{#4}\fontshape{#5}%
  \selectfont}%
\fi\endgroup%
\begin{picture}(3704,851)(1749,-81)
\put(2380,-68){\makebox(0,0)[lb]{\smash{\SetFigFont{12}{14.4}{\familydefault}{\mddefault}{\updefault}\special{ps: gsave 0 0 0 setrgbcolor}$x_1$\special{ps: grestore}}}}
\put(3140,-70){\makebox(0,0)[lb]{\smash{\SetFigFont{12}{14.4}{\familydefault}{\mddefault}{\updefault}\special{ps: gsave 0 0 0 setrgbcolor}$x_2$\special{ps: grestore}}}}
\put(4328,-81){\makebox(0,0)[lb]{\smash{\SetFigFont{12}{14.4}{\familydefault}{\mddefault}{\updefault}\special{ps: gsave 0 0 0 setrgbcolor}$x_4$\special{ps: grestore}}}}
\put(3592,-80){\makebox(0,0)[lb]{\smash{\SetFigFont{12}{14.4}{\familydefault}{\mddefault}{\updefault}\special{ps: gsave 0 0 0 setrgbcolor}$x_3$\special{ps: grestore}}}}
\end{picture}

%% file: Type2CDStratum.pstex_t
\begin{picture}(0,0)%
\includegraphics{Type2CDStratum.pstex}%
\end{picture}%
\setlength{\unitlength}{4144sp}%
\begingroup\makeatletter\ifx\SetFigFont\undefined%
\gdef\SetFigFont#1#2#3#4#5{%
  \reset@font\fontsize{#1}{#2pt}%
  \fontfamily{#3}\fontseries{#4}\fontshape{#5}%
  \selectfont}%
\fi\endgroup%
\begin{picture}(3305,833)(1654,-520)
\put(1850,-520){\makebox(0,0)[lb]{\smash{\SetFigFont{12}{14.4}{\familydefault}{\mddefault}{\updefault}\special{ps: gsave 0 0 0 setrgbcolor}$x_1$\special{ps: grestore}}}}
\put(2818,-512){\makebox(0,0)[lb]{\smash{\SetFigFont{12}{14.4}{\familydefault}{\mddefault}{\updefault}\special{ps: gsave 0 0 0 setrgbcolor}$x_2=x_3$\special{ps: grestore}}}}
\put(4197,-512){\makebox(0,0)[lb]{\smash{\SetFigFont{12}{14.4}{\familydefault}{\mddefault}{\updefault}\special{ps: gsave 0 0 0 setrgbcolor}$x_4$\special{ps: grestore}}}}
\end{picture}

%% file: Type2TD.pstex_t
\begin{picture}(0,0)%
\includegraphics{Type2TD.pstex}%
\end{picture}%
\setlength{\unitlength}{4144sp}%
\begingroup\makeatletter\ifx\SetFigFont\undefined%
\gdef\SetFigFont#1#2#3#4#5{%
  \reset@font\fontsize{#1}{#2pt}%
  \fontfamily{#3}\fontseries{#4}\fontshape{#5}%
  \selectfont}%
\fi\endgroup%
\begin{picture}(2076,1380)(1420,-610)
\put(1587,-610){\makebox(0,0)[lb]{\smash{\SetFigFont{12}{14.4}{\familydefault}{\mddefault}{\updefault}\special{ps: gsave 0 0 0 setrgbcolor}$x_1$\special{ps: grestore}}}}
\put(2217,-602){\makebox(0,0)[lb]{\smash{\SetFigFont{12}{14.4}{\familydefault}{\mddefault}{\updefault}\special{ps: gsave 0 0 0 setrgbcolor}$x_2$\special{ps: grestore}}}}
\put(2877,-595){\makebox(0,0)[lb]{\smash{\SetFigFont{12}{14.4}{\familydefault}{\mddefault}{\updefault}\special{ps: gsave 0 0 0 setrgbcolor}$x_3$\special{ps: grestore}}}}
\put(2285,575){\makebox(0,0)[lb]{\smash{\SetFigFont{12}{14.4}{\familydefault}{\mddefault}{\updefault}\special{ps: gsave 0 0 0 setrgbcolor}$x_4$\special{ps: grestore}}}}
\end{picture}

%% file: Orientation1.pstex_t
\begin{picture}(0,0)%
\includegraphics{Orientation1.pstex}%
\end{picture}%
\setlength{\unitlength}{4144sp}%
\begingroup\makeatletter\ifx\SetFigFont\undefined%
\gdef\SetFigFont#1#2#3#4#5{%
  \reset@font\fontsize{#1}{#2pt}%
  \fontfamily{#3}\fontseries{#4}\fontshape{#5}%
  \selectfont}%
\fi\endgroup%
\begin{picture}(1942,1591)(372,-943)
\put(598,-943){\makebox(0,0)[lb]{\smash{\SetFigFont{12}{14.4}{\familydefault}{\mddefault}{\updefault}\special{ps: gsave 0 0 0 setrgbcolor}$i_3$\special{ps: grestore}}}}
\put(1371,-943){\makebox(0,0)[lb]{\smash{\SetFigFont{12}{14.4}{\familydefault}{\mddefault}{\updefault}\special{ps: gsave 0 0 0 setrgbcolor}$i_4$\special{ps: grestore}}}}
\put(1138,132){\makebox(0,0)[lb]{\smash{\SetFigFont{12}{14.4}{\familydefault}{\mddefault}{\updefault}\special{ps: gsave 0 0 0 setrgbcolor}$i_1$\special{ps: grestore}}}}
\put(1148,-306){\makebox(0,0)[lb]{\smash{\SetFigFont{12}{14.4}{\familydefault}{\mddefault}{\updefault}\special{ps: gsave 0 0 0 setrgbcolor}$i_2$\special{ps: grestore}}}}
\end{picture}

%% file: Orientation2.pstex_t
\begin{picture}(0,0)%
\includegraphics{Orientation2.pstex}%
\end{picture}%
\setlength{\unitlength}{4144sp}%
\begingroup\makeatletter\ifx\SetFigFont\undefined%
\gdef\SetFigFont#1#2#3#4#5{%
  \reset@font\fontsize{#1}{#2pt}%
  \fontfamily{#3}\fontseries{#4}\fontshape{#5}%
  \selectfont}%
\fi\endgroup%
\begin{picture}(1191,1527)(547,-680)
\put(547,-680){\makebox(0,0)[lb]{\smash{\SetFigFont{12}{14.4}{\familydefault}{\mddefault}{\updefault}\special{ps: gsave 0 0 0 setrgbcolor}$i_2$\special{ps: grestore}}}}
\put(1626,-673){\makebox(0,0)[lb]{\smash{\SetFigFont{12}{14.4}{\familydefault}{\mddefault}{\updefault}\special{ps: gsave 0 0 0 setrgbcolor}$i_3$\special{ps: grestore}}}}
\put(1229,133){\makebox(0,0)[lb]{\smash{\SetFigFont{12}{14.4}{\familydefault}{\mddefault}{\updefault}\special{ps: gsave 0 0 0 setrgbcolor}$i_1$\special{ps: grestore}}}}
\end{picture}

%% file: RemainingPrincipal2.pstex_t
\begin{picture}(0,0)%
\includegraphics{RemainingPrincipal2.pstex}%
\end{picture}%
\setlength{\unitlength}{4144sp}%
\begingroup\makeatletter\ifx\SetFigFont\undefined%
\gdef\SetFigFont#1#2#3#4#5{%
  \reset@font\fontsize{#1}{#2pt}%
  \fontfamily{#3}\fontseries{#4}\fontshape{#5}%
  \selectfont}%
\fi\endgroup%
\begin{picture}(5813,965)(121,-2285)
\put(726,-2240){\makebox(0,0)[lb]{\smash{{\SetFigFont{10}{12.0}{\familydefault}{\mddefault}{\updefault}{\color[rgb]{0,0,0}$(a)$}%
}}}}
\put(2248,-2234){\makebox(0,0)[lb]{\smash{{\SetFigFont{10}{12.0}{\familydefault}{\mddefault}{\updefault}{\color[rgb]{0,0,0}$(b)$}%
}}}}
\put(3759,-2241){\makebox(0,0)[lb]{\smash{{\SetFigFont{10}{12.0}{\familydefault}{\mddefault}{\updefault}{\color[rgb]{0,0,0}$(c)$}%
}}}}
\put(5181,-2235){\makebox(0,0)[lb]{\smash{{\SetFigFont{10}{12.0}{\familydefault}{\mddefault}{\updefault}{\color[rgb]{0,0,0}$(d)$}%
}}}}
\end{picture}%

%% file: skeinrelation.pstex_t
\begin{picture}(0,0)%
\includegraphics{skeinrelation.pstex}%
\end{picture}%
\setlength{\unitlength}{4144sp}%
\begingroup\makeatletter\ifx\SetFigFont\undefined%
\gdef\SetFigFont#1#2#3#4#5{%
  \reset@font\fontsize{#1}{#2pt}%
  \fontfamily{#3}\fontseries{#4}\fontshape{#5}%
  \selectfont}%
\fi\endgroup%
\begin{picture}(4320,702)(946,93)
\put(3736,389){\makebox(0,0)[lb]{\smash{\SetFigFont{12}{14.4}{\familydefault}{\mddefault}{\updefault}\special{ps: gsave 0 0 0 setrgbcolor}$\Big)-V\Big($\special{ps: grestore}}}}
\put(2161,389){\makebox(0,0)[lb]{\smash{\SetFigFont{12}{14.4}{\familydefault}{\mddefault}{\updefault}\special{ps: gsave 0 0 0 setrgbcolor}$\Big)=V\Big($\special{ps: grestore}}}}
\put(5266,389){\makebox(0,0)[lb]{\smash{\SetFigFont{12}{14.4}{\familydefault}{\mddefault}{\updefault}\special{ps: gsave 0 0 0 setrgbcolor}$\Big)$\special{ps: grestore}}}}
\put(946,389){\makebox(0,0)[lb]{\smash{\SetFigFont{12}{14.4}{\familydefault}{\mddefault}{\updefault}\special{ps: gsave 0 0 0 setrgbcolor}$V\Big($\special{ps: grestore}}}}
\end{picture}

%% file: 4T.pstex_t
\begin{picture}(0,0)%
\includegraphics{4T.pstex}%
\end{picture}%
\setlength{\unitlength}{4144sp}%
\begingroup\makeatletter\ifx\SetFigFont\undefined%
\gdef\SetFigFont#1#2#3#4#5{%
  \reset@font\fontsize{#1}{#2pt}%
  \fontfamily{#3}\fontseries{#4}\fontshape{#5}%
  \selectfont}%
\fi\endgroup%
\begin{picture}(5109,1295)(-68,-831)
\put(128,156){\makebox(0,0)[lb]{\smash{\SetFigFont{12}{14.4}{\familydefault}{\mddefault}{\updefault}\special{ps: gsave 0 0 0 setrgbcolor}$f\Big($\special{ps: grestore}}}}
\put(2393,156){\makebox(0,0)[lb]{\smash{\SetFigFont{12}{14.4}{\familydefault}{\mddefault}{\updefault}\special{ps: gsave 0 0 0 setrgbcolor}$\Big)-f\Big($\special{ps: grestore}}}}
\put(2390,-702){\makebox(0,0)[lb]{\smash{\SetFigFont{12}{14.4}{\familydefault}{\mddefault}{\updefault}\special{ps: gsave 0 0 0 setrgbcolor}$\Big)-f\Big($\special{ps: grestore}}}}
\put(5041,154){\makebox(0,0)[lb]{\smash{\SetFigFont{12}{14.4}{\familydefault}{\mddefault}{\updefault}\special{ps: gsave 0 0 0 setrgbcolor}$\Big)$\special{ps: grestore}}}}
\put(5034,-700){\makebox(0,0)[lb]{\smash{\SetFigFont{12}{14.4}{\familydefault}{\mddefault}{\updefault}\special{ps: gsave 0 0 0 setrgbcolor}$\Big)$\special{ps: grestore}}}}
\put(-68,-706){\makebox(0,0)[lb]{\smash{\SetFigFont{12}{14.4}{\familydefault}{\mddefault}{\updefault}\special{ps: gsave 0 0 0 setrgbcolor}$=f\Big($\special{ps: grestore}}}}
\end{picture}

%% file: Universality3.pstex_t
\begin{picture}(0,0)%
\includegraphics{Universality3.pstex}%
\end{picture}%
\setlength{\unitlength}{4144sp}%
\begingroup\makeatletter\ifx\SetFigFont\undefined%
\gdef\SetFigFont#1#2#3#4#5{%
  \reset@font\fontsize{#1}{#2pt}%
  \fontfamily{#3}\fontseries{#4}\fontshape{#5}%
  \selectfont}%
\fi\endgroup%
\begin{picture}(6104,1278)(1610,-105)
\put(1914,-36){\makebox(0,0)[lb]{\smash{{\SetFigFont{10}{12.0}{\familydefault}{\mddefault}{\updefault}{\color[rgb]{0,0,0}$B_i$}%
}}}}
\put(3325,-43){\makebox(0,0)[lb]{\smash{{\SetFigFont{10}{12.0}{\familydefault}{\mddefault}{\updefault}{\color[rgb]{0,0,0}$B_j$}%
}}}}
\put(2019,694){\makebox(0,0)[lb]{\smash{{\SetFigFont{10}{12.0}{\familydefault}{\mddefault}{\updefault}{\color[rgb]{0,0,0}$p_k$}%
}}}}
\put(3344,696){\makebox(0,0)[lb]{\smash{{\SetFigFont{10}{12.0}{\familydefault}{\mddefault}{\updefault}{\color[rgb]{0,0,0}$p_{k'}$}%
}}}}
\put(5195,-50){\makebox(0,0)[lb]{\smash{{\SetFigFont{10}{12.0}{\familydefault}{\mddefault}{\updefault}{\color[rgb]{0,0,0}$B_i$}%
}}}}
\put(6606,-57){\makebox(0,0)[lb]{\smash{{\SetFigFont{10}{12.0}{\familydefault}{\mddefault}{\updefault}{\color[rgb]{0,0,0}$B_j$}%
}}}}
\put(5300,680){\makebox(0,0)[lb]{\smash{{\SetFigFont{10}{12.0}{\familydefault}{\mddefault}{\updefault}{\color[rgb]{0,0,0}$p_k$}%
}}}}
\put(7113,916){\makebox(0,0)[lb]{\smash{{\SetFigFont{10}{12.0}{\familydefault}{\mddefault}{\updefault}{\color[rgb]{0,0,0}$p_{k'}$}%
}}}}
\put(7599,1160){\makebox(0,0)[lb]{\smash{{\SetFigFont{10}{12.0}{\familydefault}{\mddefault}{\updefault}{\color[rgb]{0,0,0}...}%
}}}}
\put(4435,375){\makebox(0,0)[lb]{\smash{{\SetFigFont{10}{12.0}{\familydefault}{\mddefault}{\updefault}{\color[rgb]{0,0,0}or}%
}}}}
\end{picture}%

%% file: Universality4.pstex_t
\begin{picture}(0,0)%
\includegraphics{Universality4.pstex}%
\end{picture}%
\setlength{\unitlength}{4144sp}%
\begingroup\makeatletter\ifx\SetFigFont\undefined%
\gdef\SetFigFont#1#2#3#4#5{%
  \reset@font\fontsize{#1}{#2pt}%
  \fontfamily{#3}\fontseries{#4}\fontshape{#5}%
  \selectfont}%
\fi\endgroup%
\begin{picture}(2888,1201)(1696,-1986)
\put(4337,-1783){\makebox(0,0)[lb]{\smash{{\SetFigFont{11}{13.2}{\familydefault}{\mddefault}{\updefault}{\color[rgb]{0,0,0}...}%
}}}}
\put(4337,-1783){\makebox(0,0)[lb]{\smash{{\SetFigFont{11}{13.2}{\familydefault}{\mddefault}{\updefault}{\color[rgb]{0,0,0}...}%
}}}}
\put(3936,-1917){\makebox(0,0)[lb]{\smash{{\SetFigFont{11}{13.2}{\familydefault}{\mddefault}{\updefault}{\color[rgb]{0,0,0}$p_{k}'$}%
}}}}
\put(4210,-1552){\makebox(0,0)[lb]{\smash{{\SetFigFont{11}{13.2}{\familydefault}{\mddefault}{\updefault}{\color[rgb]{0,0,0}$B_i$}%
}}}}
\put(3915,-941){\makebox(0,0)[lb]{\smash{{\SetFigFont{11}{13.2}{\familydefault}{\mddefault}{\updefault}{\color[rgb]{0,0,0}$p_k$}%
}}}}
\put(2492,-1803){\makebox(0,0)[lb]{\smash{{\SetFigFont{11}{13.2}{\familydefault}{\mddefault}{\updefault}{\color[rgb]{0,0,0}...}%
}}}}
\put(2492,-1803){\makebox(0,0)[lb]{\smash{{\SetFigFont{11}{13.2}{\familydefault}{\mddefault}{\updefault}{\color[rgb]{0,0,0}...}%
}}}}
\put(2091,-1938){\makebox(0,0)[lb]{\smash{{\SetFigFont{11}{13.2}{\familydefault}{\mddefault}{\updefault}{\color[rgb]{0,0,0}$p_{k}'$}%
}}}}
\put(2364,-1572){\makebox(0,0)[lb]{\smash{{\SetFigFont{11}{13.2}{\familydefault}{\mddefault}{\updefault}{\color[rgb]{0,0,0}$B_i$}%
}}}}
\put(2254,-1198){\makebox(0,0)[lb]{\smash{{\SetFigFont{11}{13.2}{\familydefault}{\mddefault}{\updefault}{\color[rgb]{0,0,0}$p_k$}%
}}}}
\end{picture}%

%% file: Universality8.pstex_t
\begin{picture}(0,0)%
\includegraphics{Universality8.pstex}%
\end{picture}%
\setlength{\unitlength}{4144sp}%
\begingroup\makeatletter\ifx\SetFigFont\undefined%
\gdef\SetFigFont#1#2#3#4#5{%
  \reset@font\fontsize{#1}{#2pt}%
  \fontfamily{#3}\fontseries{#4}\fontshape{#5}%
  \selectfont}%
\fi\endgroup%
\begin{picture}(4798,1327)(-325,-433)
\put(1931,155){\makebox(0,0)[lb]{\smash{{\SetFigFont{10}{12.0}{\familydefault}{\mddefault}{\updefault}{\color[rgb]{0,0,0}or}%
}}}}
\put(-13,-381){\makebox(0,0)[lb]{\smash{{\SetFigFont{10}{12.0}{\familydefault}{\mddefault}{\updefault}{\color[rgb]{0,0,0}$B_k$}%
}}}}
\put(1259,-372){\makebox(0,0)[lb]{\smash{{\SetFigFont{10}{12.0}{\familydefault}{\mddefault}{\updefault}{\color[rgb]{0,0,0}$B_{k'}$}%
}}}}
\put(574,762){\makebox(0,0)[lb]{\smash{{\SetFigFont{10}{12.0}{\familydefault}{\mddefault}{\updefault}{\color[rgb]{0,0,0}$p_i$}%
}}}}
\put(2576,-384){\makebox(0,0)[lb]{\smash{{\SetFigFont{10}{12.0}{\familydefault}{\mddefault}{\updefault}{\color[rgb]{0,0,0}$B_k$}%
}}}}
\put(3849,-375){\makebox(0,0)[lb]{\smash{{\SetFigFont{10}{12.0}{\familydefault}{\mddefault}{\updefault}{\color[rgb]{0,0,0}$B_{k'}$}%
}}}}
\put(3131,759){\makebox(0,0)[lb]{\smash{{\SetFigFont{10}{12.0}{\familydefault}{\mddefault}{\updefault}{\color[rgb]{0,0,0}$p_i$}%
}}}}
\end{picture}%

%% file: Universality9.pstex_t
\begin{picture}(0,0)%
\includegraphics{Universality9.pstex}%
\end{picture}%
\setlength{\unitlength}{4144sp}%
\begingroup\makeatletter\ifx\SetFigFont\undefined%
\gdef\SetFigFont#1#2#3#4#5{%
  \reset@font\fontsize{#1}{#2pt}%
  \fontfamily{#3}\fontseries{#4}\fontshape{#5}%
  \selectfont}%
\fi\endgroup%
\begin{picture}(2292,1524)(-557,-569)
\put(-68,-342){\makebox(0,0)[lb]{\smash{{\SetFigFont{10}{12.0}{\familydefault}{\mddefault}{\updefault}{\color[rgb]{0,0,0}$B_{k}$}%
}}}}
\put(278,-519){\makebox(0,0)[lb]{\smash{{\SetFigFont{10}{12.0}{\familydefault}{\mddefault}{\updefault}{\color[rgb]{0,0,0}$\widetilde{B}_{k}$}%
}}}}
\put(1478,-523){\makebox(0,0)[lb]{\smash{{\SetFigFont{10}{12.0}{\familydefault}{\mddefault}{\updefault}{\color[rgb]{0,0,0}$\widetilde{B}_{k'}$}%
}}}}
\put(1132,-346){\makebox(0,0)[lb]{\smash{{\SetFigFont{10}{12.0}{\familydefault}{\mddefault}{\updefault}{\color[rgb]{0,0,0}$B_{k'}$}%
}}}}
\put(472,831){\makebox(0,0)[lb]{\smash{{\SetFigFont{10}{12.0}{\familydefault}{\mddefault}{\updefault}{\color[rgb]{0,0,0}$p_i$}%
}}}}
\end{picture}%